\documentclass{amsart}
\usepackage[british,ngerman]{babel}
\usepackage{amsmath,amssymb}
\usepackage[pagebackref]{hyperref}

\newcommand{\N}{\mathbb N}
\newcommand{\Z}{\mathbb Z}
\newcommand{\Q}{\mathbb Q}
\newcommand{\C}{\mathbb C}
\newcommand{\fa}{\mathfrak a}
\newcommand{\fb}{\mathfrak b}
\newcommand{\fc}{\mathfrak c}

\newcommand{\ggT}{\operatorname{ggT}}

\newcommand{\grad}{\operatorname{grad}}

\newcommand{\impl}{\Longrightarrow}

\newcommand{\ct}{\operatorname{ct}}

\def\proofname{Beweis}
\newtheorem{thm}{Satz}
\newtheorem{prop}[thm]{Proposition}
\newtheorem{lem}[thm]{Lemma}

\newenvironment{bew}{\begin{proof}}{\end{proof}}

\renewcommand{\theHsection}{de.\arabic{section}}


\begin{document}

\title{Zur Zahlentheorie der Griechen}

\author{Franz Lemmermeyer}

\begin{abstract}
Im zweiten Teil (vgl. \cite{LSb}) unseres Artikels \"uber den
Euklidischen Vierzahlensatz gehen wir vom Begriff der
Inkommensurabilit\"at aus und leuchten die Irrationalit\"atsbeweise
algebraischer Zahlen aus. Dabei treffen wir auf das Gau\ss{}sche
Lemma, die ganze Abgeschlossenheit und Dedekinds Prager Satz, um
schlie\ss{}lich wieder beim Vierzahlensatz zu landen.
\end{abstract}

\maketitle

\begin{center}
{\Large Teil II. Gau\ss{}sche Lemmata und Rieszsche Ringe}
\end{center}
\medskip

\setcounter{section}{7}

\section{Noch einmal die Alten Griechen: Inkommensurabilit\"at}

Bekanntlich wird den Pythagoreern auch die Entdeckung der
Inkommensurabilit\"at zugeschrieben, also der Beobachtung, dass
z.B. Diagonale $d$ und Seite $s$ eines Quadrats kein rationales Verh\"altnis
zueinander haben, dass sich also nicht $d:s = p:q$ f\"ur nat\"urliche
Zahlen $p, q$ schreiben l\"asst. F\"ur diese Behauptung gibt es eine
Unzahl von Beweisen, von denen zwei bereits auf die Griechen
zur\"uckgehen:
\begin{itemize}
\item Der eine Beweis folgert aus $\sqrt{2} = \frac pq$, dass $p$ und
      $q$ beide gerade sein m\"ussen, was der K\"urzbarkeit von
      Br\"uchen (bzw. Proportionen) widerspricht. Dieser Beweis geht
      auf Aristoteles zur\"uck und ist in manchen (\"alteren) Ausgaben
      der Elemente als Proposition X.117 enthalten; er gilt heute aber
      als sp\"ater eingef\"ugter Zusatz.
\item Der andere Beweis (vgl. Hasse \cite{Has} und Apostol \cite{Apo})
      ist geometrischer Natur und basiert auf der ``Wechselwegnahme'',
      dem Vater des Euklidischen Algorithmus.
\end{itemize}
Dass bereits die Pythagoreer die Existenz inkommensurabler Strecken
entdeckten, gilt als sicher; dass Hippasus diese Entdeckung gemacht
(bzw. bekannt gemacht) hat, ist dagegen ungewiss, ebenso die
Behauptung, er sei durch das Studium des Pentagramms darauf
gef\"uhrt worden (sh. von Fritz \cite{vFr}). Die hin und wieder
anzutreffende Behauptung, dies h\"atte eine ``Grundlagenkrise'' der
griechischen Mathematik verursacht, ist schlie{\ss}lich durch keinerlei
historische Fakten mehr zu belegen (muss deswegen aber nicht falsch sein).

Die Inkommensurabilit\"at von Strecken hat sich im Zuge der Erweiterung
des Zahlbegriffs von einem geometrischen zu einem arithmetischen Problem
gewandelt: das Verh\"altnis zweier Strecken ist f\"ur uns einfach eine
Zahl, und die Inkommensurabilit\"at von Seite und Diagonale eines
Quadrats die Irrationalit\"at der Quadratwurzel aus $2$.

Wenn man wei{\ss}, dass der Ring $\Z$ faktoriell ist, wird der Beweis
der Irrationalit\"at von $\sqrt{2}$ trivial: aus $2q^2 = p^2$ liest
man ab, dass $2$ die linke Seite in ungerader, die rechte in gerader
Potenz teilt: Widerspruch.

Derselbe Beweis zeigt allgemeiner\footnote{Ohne eindeutige Zerlegung
geht das schief: im Monoid $M$ aller
nat\"urlichen Zahlen $\equiv 0, 1, 2 \bmod 4$ ist $9$ das Quadrat des
Bruchs $\frac62$ im Quotientenmonoid, aber kein Quadrat einer Zahl
in $M$. Entsprechend ist $2i$ kein Quadrat in $\Z[2i] = \Z \oplus 2i\Z$,
wohl aber eines im Quotientenk{\"o}rper $\Q(i)$ wegen
$2i = (1+i)^2 = (\frac{2+2i}{2})^2$.}:

\begin{thm}\label{Sm}
F\"ur $m \in \Z$ ist $\sqrt{m}$ genau dann rational, wenn $m$ das
Quadrat einer ganzen Zahl ist.
\end{thm}
Auch dies ist ein Satz, der im Wesentlichen bereits den Griechen
bekannt war; allerdings ist er so explizit nicht in Euklids Elementen
zu finden.

Anstatt den Beweis der Irrationalit\"at von $\sqrt{2}$ mit Hilfe der
eindeutigen Primfaktorzerlegung zu f\"uhren, k{\"o}nnen wir ihn auch
durch Ausnutzen der euklidischen Division mit Rest erhalten (was
angesichts der Tatsache, dass Euklidische Ringe faktoriell sind,
nat\"urlich keinesfalls verwundern sollte):
Sei $m \in \N$ kein Quadrat; dann gibt es ein $n \in \N_0$ mit
$n^2 < m < (n+1)^2$, also mit $0 < \sqrt{m}-n < 1$. Ist $\sqrt{m}$
rational, so k{\"o}nnen wir $\sqrt{m} = \frac pq$ mit minimalem
$q \in \N$ schreiben. Dann ist $r := q(\sqrt{m}-n) = p - nq$ ebenfalls
ganz, und wegen $0 < \sqrt{m}-n < 1$ ist $0 < r < q$. Nun ist
$r\sqrt{m} =  q(\sqrt{m}-n)\sqrt{m} = mq-nq\sqrt{m} = mq - np$ ganz;
wegen $\sqrt{m} = \frac{r\sqrt{m}}{r}$ widerspricht dies aber der
Minimalit\"at von $q$. F\"ur Verallgemeinerungen dieser
Schlussweise\footnote{Dieser
Beweis geht im Wesentlichen auf Dedekind \cite{DedSIZ} zur\"uck;
f\"ur eine geometrische Einkleidung des Beweises durch Haj\'os sh.
\cite{Sur2}. Varianten dieses Beweises finden sich u.a. in
\cite{Blo,Fine,Sag,Vau,Wat}.} sowie andere Irrationalit\"atsbeweise
sh. z.B. die \"Ubungen in \cite[Kap. 1]{LL}.

Hier ist ein Beweis, welcher auf der Tatsache basiert, dass $\Z$
ein Bezout-Ring ist; er zeigt gleich allgemeiner, dass eine rationale
$n$-te Wurzel einer ganzen Zahl $m$ bereits ganz ist: sei
$\sqrt[n]{m} = \frac pq$ mit teilerfremden ganzen Zahlen $p$ und $q$,
also $(p,q) = (1)$. Dann ist wegen $p^n = mq^n$
$$ (1) = (p,q)^n \subseteq (p^n,q) = (mq^n,q) \subseteq (q), $$
folglich $(q) = (1)$, also $q = \pm 1$ und $\sqrt[n]{m} = \pm p$
eine ganze Zahl.

So h\"ubsch solche Beweise auch sind, verschleiern sie in gewisser
Weise doch den eigentlichen Kern der Angelegenheit; diesen wollen wir
im n\"achsten Abschnitt offenlegen.

\section{Ganze Abgeschlossenheit}

Der Begriff der ganzen Abgeschlossenheit begegnet Studierenden in der
Regel in ein\-f\"uh\-ren\-den Algebravorlesungen: ist $R$ ein Integrit\"atsring
mit Quotientenk{\"o}rper $K$, so nennt man $t \in K$ ganz \"uber $R$,
wenn $t$ Nullstelle eines normierten Polynoms $f \in R[X]$ ist.
Ist jedes $t \in K$, das ganz \"uber $R$ ist, bereits in $R$ enthalten,
so nennt man $R$ ganz abgeschlossen (in seinem Quotientenk{\"o}rper).

Sind $t_1, t_2$ (in einem Erweiterungsk\"orper von $K$) ganz \"uber $R$,
so sind auch $t_1 + t_2$ und $t_1t_2$ ganz \"uber $R$. Der klassische
Beweis \"uber Resultanten wurde im Laufe der Zeit ersetzt durch Beweise
via endlich erzeugter Moduln; im Gewande der Tensorsprache kommt der
klassische Beweis aber an Eleganz fast an den Modulbeweis heran: $t_j$
ist als Nullstelle eines normierten Polynoms vom Grad $m_j$ mit
Koeffizienten aus $R$ auch Eigenwert einer $m_j \times m_j$-Matrix $T_j$
(der Begleitmatrix) mit Eintr\"agen aus $R$. Man zeigt nun leicht, dass
$t_1t_2$ Eigenwert des Kroneckerprodukts $T_1 \otimes T_2$ und damit
(als Nullstelle von dessen charakteristischem Polynom) ganz \"uber $R$
ist. Entsprechend folgt, dass $t_1 + t_2$ Eigenwert von
$T_1 \otimes I_{m_2} + I_{m_1} \otimes T_2$ ist, wo die $I_k$
Einheitsmatrizen der Dimension $k$ sind (sh. \cite{KMS}).

Der Satz, dass faktorielle Ringe ganz abgeschlossen sind, geh{\"o}rt zum
Grundwissen von Algebrastudenten; der \"ubliche Beweis liefert sogar
die st\"arkere Aussage, dass GGT-Ringe ganz abgeschlossen sind. Aber
selbst in Vorlesungen \"uber algebraische Zahlentheorie wird selten
klar herausgearbeitet, wie grundlegend der Begriff der ganzen
Abgeschlossenheit f\"ur die Entwicklung dieser Disziplin gewesen
ist, und welche Bedeutung er f\"ur die moderne Algebra gewonnen hat,
seitdem Dedekind und Emmy Noether
ihn\footnote{Kummer hat den Begriff der ganzen algebraischen Zahl
nicht gekannt, und die Bedeutung der ganzen Abgeschlossenheit ist
ihm nicht klar geworden. Seine Konstruktion der idealen Zahlen in
Kreisk{\"o}rpern hat deswegen funktioniert, weil die Ringe, die
Kummer betrachtete, ganz abgeschlossen waren; allerdings hat sein
Aufbau seiner Theorie der idealen Zahlen eine L\"ucke, die erst
Dedekind geschlossen hat -- unter Ausnutzung der ganzen
Abgeschlossenheit (vgl. \cite{LJ}).

Auch Dirichlet und Eisenstein blieb die ganze Abgeschlossenheit
verborgen: Dirichlet bewies seinen Einheitensatz in Ringen der
Form $\Z[\alpha]$, wo $\alpha$ eine ganze algebraische Zahl ist,
und aus einer Bemerkung Eisensteins folgt, dass die ganzen
algebraischen Zahlen einen Ring bilden. Heine \cite{Heine} definierte
ganze algebraische Zahlen als solche, die sich mittels Ringoperationen
und wiederholten Ziehens von $n$-ten Wurzeln aus den ganzen Zahlen
gewinnen lassen.

Dedekind zeigte, dass in den Ringen ganzer Zahlen eines Zahlk\"orpers
$K \subseteq \C$ jedes Ideal $\ne (0)$ eindeutig Produkt von
Primidealen ist; schlie\ss{}lich gelang es Emmy Noether 1926 zu zeigen,
dass dieser Satz in allen Integrit\"atsringen gilt, in denen gewisse
Axiome erf\"ullt sind -- eines davon ist die ganze Abgeschlossenheit.}
aus dem Aufbau der Dedekindschen Idealtheorie herauskristallisiert
haben. In der Tat zeigt ein Blick auf die folgende Abbildung, dass
die meisten der in den einf\"uhrenden Algebravorlesungen auftauchenden
Klassen von Ringen\footnote{Bisher noch nicht definierte Begriffe
wie Schreier oder $P^4$ werden in den Abschnitten \ref{Ab12} und
\ref{Ab13} erkl\"art.} ganz abgeschlossen sind:

\begin{figure}[h!]
$$ \begin{array}{cccccccccccc}
\text{Euklid} &         &       &    & & & \text{Dedekind} & & \\
    \cap      &         &       &    & & & \cap & & \\
 \text{HI}    & \subset & \text{Bezout} & & \subset & & \text{Pr\"ufer} & & \\
    \cap      &         &   \cap     &         & & & \cap & & \\
\text{ZPE} & \subset & \text{GGT} & \subset & \text{Schreier}
              & \subset & \text{ganz abg.} &  \\
    \cap      &         &      &       &  \cap   &  &  & & \\
\text{atomar} &         &      &   & \text{Riesz} & \subset & P^4 &
    \subseteq & \text{AP} \\
\end{array}$$
\caption{Klassen von Ringen}\label{FR}
\end{figure}

Die ganze Abgeschlossenheit von $\Z$ ist auch der eigentliche
Grund f\"ur die G\"ultigkeit von Satz \ref{Sm}. Ist n\"amlich
$n \in \Z$ das Quadrat einer rationalen Zahl, also z.B.
$\sqrt n = \frac pq$, so ist $\frac pq$ Nullstelle des quadratischen
Polynoms $X^2 - n \in \Z[X]$. Da $\Z$ ganz abgeschlossen ist, muss
$\frac pq \in \Z$ sein.

Es ist daher keineswegs erstaunlich, dass sich die Irrationalit\"at
von $\sqrt{2}$ durch Ausnutzen der Faktorialit\"at von $\Z$ nachweisen
l\"asst. Weiter kann man versuchen, derartige Irrationalit\"atsbeweise
mit all denjenigen S\"atzen zu f\"uhren, aus denen ebenfalls die ganze
Abgeschlossenheit von $\Z$ folgt. Ein Beispiel hierf\"ur ist der
folgende Beweis:

\begin{align*}
 \sqrt{m} \text{ ist rational}
     & \iff X^2 - m \text{ ist reduzibel in } \Q[X] \\
     & \iff X^2 - m \text{ ist reduzibel in } \Z[X] \\
     & \iff X^2 - m = (X-n)(X+n) \text{ f\"ur ein } n \in \Z \\
     & \iff m = n^2 \text{ f\"ur ein } n \in \Z.
\end{align*}
Der \"Ubergang von $\Q[X]$ zu $\Z[X]$ geschieht dabei mit
dem Gau{\ss}schen Hebungslemma\footnote{Der Name soll andeuten,
dass sich Faktorisierungen von Polynomen \"uber $\Q$ zu
Faktorisierungen \"uber $\Z$ hochheben lassen.}
f\"ur normierte Polynome:

\begin{lem}\label{GHL1}
Sind $g, h \in \Q[X]$ normierte Polynome mit $gh = f \in \Z[X]$,
dann sind $g, h \in \Z[X]$.
\end{lem}

In dieser Form taucht das Gau{\ss}sche Hebungslemma in den Disquisitiones
Arithmeticae auf (n\"amlich in \cite[Art. 42]{DA}), und wird dort zum
Nachweis der Irreduzibilit\"at der Kreisteilungspolynome
$F_p(X) = \frac{X^p-1}{X-1}$ f\"ur prime $p$ verwendet.
Zur Abgrenzung von anderen Ergebnissen, die ebenfalls unter dem
Namen ``Gau{\ss}sches Lemma'' bekannt sind\footnote{Es gibt mindestens
drei solcher Resultate: das oben angef\"uhrte Hebungslemma, und je eines
aus der Theorie der quadratischen Reste bzw. der quadratischen Formen.
Alle drei finden sich in den Disquisitiones Arithmeticae.

Weiter wird das Gau{\ss}sche Hebungslemma in $\Z$ oft anders formuliert,
und diese Versionen sind zwar \"uber $\Z$, nicht aber \"uber allgemeinen
Integrit\"atsringen zu dem oben formulierten Hebungslemma \"aquivalent.},
werden wir Lemma \ref{GHL1} das Gau{\ss}sche Hebungslemma (f\"ur
normierte Polynome) nennen.

Mit Lemma \ref{GHL1} l\"asst sich viel mehr zeigen als nur die
Irrationalit\"at gewisser Quadratwurzeln: ist n\"amlich $\alpha$
Nullstelle eines normierten Polynoms mit ganzen Koeffizienten (also
eine ganze algebraische Zahl), so ist $\alpha$ eine ganze Zahl oder
irrational. In der Tat: ist $\alpha$ rationale Nullstelle des
normierten Polynoms $f \in \Z[X]$, dann ist $f(X) = (X-\alpha)g(X)$,
und $g \in \Q[X]$ normiert. Das Gau{\ss}sche Hebungslemma f\"ur normierte
Polynome zeigt dann $X - \alpha \in \Z[X]$, also $\alpha \in \Z$.

Das Gau{\ss}sche Hebungslemma f\"ur normierte Polynome ist ein
Spezialfall von

\begin{thm}[Gau{\ss}sches Hebungslemma]\label{LGHL}
Ist $F \in \Z[X]$ ein Polynom und gibt es Polynome $g, h \in \Q[X]$
mit $F = gh$, dann existiert ein $c \in \Q^\times$ mit der Eigenschaft,
dass $G = cg, H = \frac1c h \in \Z[X]$ ist. Jede Faktorisierung
$F = gh$ in $\Q[X]$ l\"asst sich daher zu einer Faktorisierung
$F = GH$ in $\Z[X]$ ``hochheben''.
\end{thm}

Dieser Satz ist auch in heutigen Algebravorlesungen noch ein wichtiges
Hilfsmittel beim Nachweis der Irreduzibilit\"at von Polynomen
und wird beispielsweise beim Beweis des Eisensteinschen Kriteriums
benutzt: mit seiner Hilfe zeigt man n\"amlich, dass ein primitives
Polynom vom Grad $\ge 1$ in $\Z[X]$ genau dann irreduzibel in $\Q[X]$
ist, wenn es irreduzibel in $\Z[X]$ ist.

Ist $R$ ein Integrit\"atsbereich  mit Quotientenk{\"o}rper $K$, so
sagen wir, in $R[X]$ gelte das Gau{\ss}sche Hebungslemma, wenn sich jede
Faktorisierung eines $F \in R[X]$ von $K[X]$ nach $R[X]$ hochheben
l\"asst, wenn es also zu $F = gh$ mit $g, h \in K[X]$ stets ein
$c \in K^\times$ gibt mit $cg, \frac1c h \in R[X]$. Wir stellen uns
nun die Frage, in welchen Ringen die obigen Versionen des
Gau\ss{}schen Hebungslemmas gelten. F\"ur normierte
Polynome ist diese Frage relativ leicht zu beantworten:

\begin{thm}\label{GLDP}
Sei $R$ ein Integrit\"atsring mit Quotientenk{\"o}rper $K$. Dann sind
folgende Aussagen gleichbedeutend:
\begin{enumerate}
\item $R$ ist ganz abgeschlossen.
\item In $R[X]$ gilt das Gau{\ss}sche Hebungslemma f\"ur normierte
      Polynome: \\
      Sind $g, h \in K[X]$ normierte Polynome mit $gh = f \in R[X]$,
      dann sind $g, h \in R[X]$.
\item In $R[X]$ gilt Dedekinds Prager Satz: sind \\
      \begin{equation}\label{EDPS}
       b(X) = \sum_{i=0}^m b_i X^i \quad \text{und } \quad
       c(X) = \sum_{j=0}^n c_jX^j
      \end{equation}
      Polynome in $R[X]$, und ist $a \mid b(X)c(X)$ f\"ur ein $a \in R$,
      so gilt $a \mid b_ic_j$ f\"ur alle $0 \le i \le m$, $0 \le j \le n$.
\end{enumerate}
\end{thm}

Bevor wir den Beweis von Satz \ref{GLDP} geben und danach die
G\"ultigkeit des allgemeinen Hebungslemmas untersuchen, machen
wir noch einige Bemerkungen zur Bedeutung von Dedekinds Prager Satz.

\section{Dedekinds Prager Satz}

Historisch spielte Dedekinds Prager Satz\footnote{Diesen Satz hat
Dedekind in der Zeitschrift der Prager Mathematischen Gesellschaft
publiziert, und heute hat sich der Name ``Prager Satz'' daf\"ur
eingeb\"urgert.} eine wichtige Rolle
beim Aufbau der Dedekindschen Idealtheorie. Dedekind empfand
das Argumentieren mit Polynomen beim Aufbau seiner Idealtheorie
als unangemessen\footnote{Vgl. Dedekinds Kommentare zum Hurwitzschen
Aufbau der Idealtheorie in \cite{DedH}, so auf S. 53 (``weil die
Benutzung der Funktionen von Variablen mir immer als ein der Sache
fremdes Hilfsmittel erscheint'') und S. 55 (``Aus denselben
Gr\"unden konnte der oben erw\"ahnte Beweis [\ldots] mich noch
nicht v\"ollig befriedigen, weil durch die Einmischung der
Funktionen von Variablen die Reinheit der Theorie nach meiner
Ansicht getr\"ubt wird'').} und suchte nach anderen Wegen. Hurwitz
\cite{Hurw} ver{\"o}ffentlichte dann einen Zugang zur Idealtheorie,
der auf Dedekinds Prager Satz aufbaute (welcher Hurwitz damals
nicht bekannt war), woraufhin Dedekind erkl\"arte, dass ihm
dieser Weg bekannt sei, er aus methodischen Gr\"unden aber seinem
eigenen Aufbau den Vorzug gegeben habe. Hilberts Zahlbericht
\cite{HZb} und Landaus\footnote{Landau \cite[S. 145]{Lan} nennt
Dedekinds Prager Satz einen ``Satz von Kronecker''; bei Cohn
\cite{CAN} hei\ss{}t er ``Hurwitz's Lemma''. Walter Strobl weist
in seiner Diplomarbeit \cite[S. 157--161]{Str} darauf hin, dass
Bachmann in Bd. III und Bd. V seines f\"unfb\"andigen Lehrwerks \"uber
Zahlentheorie detailliert auf den Prager Satz eingeht. Vergleiche
auch die Einf\"uhrung in Edwards \cite{EdDiv}.} Lehrbuch \cite{Lan}
\"uber algebraische Zahlen benutzen beide den Hurwitzschen
Zugang\footnote{Dedekinds Satz dient dem Nachweis, dass es zu jedem
Ideal $\fa \ne (0)$ in einer Maximalordnung eines Zahlk{\"o}rpers
ein Ideal $\fb \ne (0)$ gibt derart, dass $\fa\fb$ ein Hauptideal ist.},
w\"ahrend man heute eher mit Lokalisierungen arbeitet (was sicher den
Beifall Dedekinds gefunden h\"atte).

Dedekind \cite{Ded} gab zwei Beweise seines Prager Satzes;
den ersten findet man bei Landau \cite{Lan} (Hilbert gibt
in seinem Zahlbericht \cite{HZb} keinen Beweis, sondern verweist
auf Arbeiten von Dedekind \cite{Ded}, Hurwitz \cite{Hurw},
Kronecker \cite{Kron} und Mertens \cite{Mert}), und
der zweite Beweis f\"uhrt auf das Lemma von Dedekind--Mertens,
dessen Geschichte und Weiterentwicklung u.a. durch Artin und
Rees einen eigenen Artikel verdient h\"atte (vgl. dazu
\cite{Coqu,CDLQ}).

Wir wollen hier einen dritten Beweis angeben, der auf der
Beobachtung beruht, dass Dedekinds Prager Satz genau
in ganz abgeschlossenen Ringen gilt\footnote{Dass Emmy Noether
diesen Satz \ref{GLDP} gekannt hat, wird aus ihren Bemerkungen
zu Dedekinds Artikel \cite{DedH} in seinen Werken klar.}; damit
ist der Prager Satz \"uber die ganze Abgeschlossenheit \"aquivalent
zum Gau\ss{}schen Hebungslemma f\"ur normierte Polynome. Der
entscheidende Schritt (1) $\impl$ (3) beruht auf einem klassischen
Lemma, das im Wesentlichen auf Kronecker zur\"uckgeht:

\begin{lem}[Kronecker]\label{LKron}
Sei $R$ ein Integrit\"atsring mit Quotientenk\"orper $K$, und sei
$f(X) = \lambda \prod_{k=1}^N (X - \theta_k) \in R[X]$ mit $\lambda \ne 0$
und $\theta_k$ in einem Erweiterungsk\"orper $L$ von $K$. Dann ist
f\"ur jede Teilmenge $I \subseteq \{1, \ldots, N\}$ das Element
$\lambda \theta_I$, wo $\theta_I = \prod_{k \in I} \theta_k$, ganz
\"uber $R$.
\end{lem}

\begin{bew}
Seien $e_1, \ldots, e_N$ die elementarsymmetrischen Funktionen der
$\theta_k$; wegen $f = \lambda X^N - \lambda e_1 X^{N-1} \pm \ldots$
ist $\lambda e_k \in R$ f\"ur alle $k$, und nat\"urlich ist
$\lambda \in R$. F\"ur festes $r = |I|$ betrachten wir das Polynom
$$ P(Y) = \prod_{|J| = r} (Y - \lambda\theta_J), $$
wobei $J$ alle $r$-elementigen Teilmengen von $\{1, \ldots, N\}$
durchl\"auft. Der Koeffizient von $Y^{M-s}$ (mit $M = \binom Nr$) ist
$(-1)^s\lambda^s E_s$, wo $E_s$ die $s$-te elementarsymmetrische
Funktion der $\theta_J$ ist. $E_s$ ist ein symmetrisches Polynom in
$\theta_1, \ldots, \theta_N$ mit ganzen Koeffizienten, in dem jede
Variable h\"ochstens in der $s$-ten Potenz vorkommt (jedes $\theta_J$
enth\"alt jedes $\theta_k$ h\"ochstens einfach). Nach dem Hauptsatz
\"uber symmetrische Funktionen (in der durch den \"ublichen Beweis
gelieferten pr\"azisierten Form) ist $E_s$ daher ein Polynom vom
Gesamtgrad $\le s$ in $e_1, \ldots, e_N$ mit ganzen Koeffizienten.
F\"ur jedes Monom $e_{k_1} \cdots e_{k_v}$ mit $v \le s$ ist aber
$\lambda^s e_{k_1} \cdots e_{k_v} = \lambda^{s-v} (\lambda e_{k_1})
\cdots (\lambda e_{k_v}) \in R$. Also ist $P$ normiert mit
Koeffizienten in $R$, und $\lambda\theta_I$ ist als Nullstelle
von $P$ ganz \"uber $R$.
\end{bew}

\begin{bew}[von Satz \ref{GLDP}]
(1) $\impl$ (2): Sei $L$ ein Zerf\"allungsk\"orper von $f$ \"uber $K$;
dort zerfallen auch $g$ und $h$ in Linearfaktoren, etwa
$g(X) = \prod_i (X - \alpha_i)$ und $h(X) = \prod_j (X - \beta_j)$.
Wegen $f(\alpha_i) = f(\beta_j) = 0$ sind die $\alpha_i$ und $\beta_j$
ganz \"uber $R$. Also sind auch die elementarsymmetrischen
Funktionen in den $\alpha_i$ bzw. $\beta_j$, d.h. die
Koeffizienten von $g$ und $h$, ganz \"uber $R$. Da diese in $K$ liegen
und $R$ ganz abgeschlossen ist, liegen sie bereits in $R$, und es folgt
$g, h \in R[X]$ (dieser Beweis stammt aus \cite{CS}).

(2) $\impl$ (1): Sei $t \in K$ ganz \"uber $R$, also Nullstelle
eines normierten Polynoms $f \in R[X]$. Dann ist $f(X) = (X-t)g(X)$
mit normiertem $g \in K[X]$, und nach (2) ist $X - t \in R[X]$, also
$t \in R$.

(1) $\impl$ (3): Sei $a \mid b(X)c(X)$. Ist $a = 0$ oder $b = 0$ oder
$c = 0$, so ist nichts zu zeigen (im Falle $a = 0$ ist $bc = 0$). Seien
also $a, b, c \ne 0$, $b_m \ne 0 \ne c_n$, und $f = \frac{bc}a \in R[X]$
mit Leitkoeffizient $\lambda = \frac{b_mc_n}a$. In einem
Zerf\"allungsk\"orper $L$ von $bc$ \"uber $K$ ist
$b = b_m\prod_{k=1}^m (X-\beta_k)$, $c = c_n \prod_{l=1}^n (X-\gamma_l)$,
also $f = \lambda \prod_{k=1}^{m+n} (X-\theta_k)$ mit
$(\theta_1, \ldots, \theta_{m+n}) = (\beta_1, \ldots, \beta_m,
\gamma_1, \ldots, \gamma_n)$. Nun ist
$$ \frac{b_ic_j}a = \lambda \cdot \frac{b_i}{b_m} \cdot \frac{c_j}{c_n}
   = \pm \lambda \, e_{m-i}(\beta_1, \ldots, \beta_m) \,
           e_{n-j}(\gamma_1, \ldots, \gamma_n), $$
und das Ausmultiplizieren der rechten Seite liefert eine Summe von
Termen der Form $\lambda \theta_I$ mit $I \subseteq \{1, \ldots, m+n\}$.
Nach Lemma \ref{LKron} ist daher $\frac{b_ic_j}a \in K$ ganz \"uber $R$,
also in $R$, d.h. $a \mid b_ic_j$.

(3) $\impl$ (1): Sei $t \in K$ ganz \"uber $R$, also Nullstelle
eines normierten Polynoms $f \in R[X]$. Dann ist $f(X) = (X-t)g(X)$
mit normiertem $g \in K[X]$; seien $a, b \in R \setminus \{0\}$ mit
$at \in R$ und $bg \in R[X]$. Also ist
$[a(X-t)][bg(X)] = abf(X)$ durch $ab$ teilbar, und Dedekinds Prager
Satz zeigt, dass $ab \mid atb$ in $R$ gilt (auf der rechten Seite
steht bis aufs Vorzeichen das Produkt des konstanten Glieds von
$a(X-t)$ und des h{\"o}chsten Koeffizienten $b$ von $bg(X)$); daraus
folgt $t = \frac{atb}{ab} \in R$.
\end{bew}

\noindent{\bf Bemerkung.} Ein direkter Schluss von (2) auf (3) nach
dem Muster ``$f = \frac{b}{a} \cdot c$ ist eine Faktorisierung in
$K[X]$, also l\"asst sie sich hochheben'' ist nicht m\"oglich, da
$b$ und $c$ im Allgemeinen nicht normiert sind: dazu br\"auchte man
das allgemeine Hebungslemma, und dieses gilt, wie wir sehen werden,
nur in Schreierschen Ringen.
\medskip

Dedekinds Prager Satz spielt in der heutigen Algebra und Zahlentheorie
praktisch keine Rolle mehr; das ist, wie Dedekind selbst bemerkt
hat, nur konsequent: er ist ja \"aquivalent zur ganzen Abgeschlossenheit
des Grundrings, und anstatt rechnerische Beweise mit dem Prager
Satz zu f\"uhren, sollte man nach Dedekind konzeptionelle Beweise
mit abstrakten Prinzipien wie der ganzen Abgeschlossenheit bevorzugen.
Die Geschichte hat ihm da, wie bereits Emmy Noether in ihren
Kommentaren zu Dedekinds Werken schrieb, Recht gegeben, auch wenn
insbesondere Hilbert in seinem Zahlbericht dem Hurwitzschen Aufbau
der Idealtheorie noch den Vorzug gegeben hat. Dieser ``Fehlgriff''
(vom Standpunkt der abstrakten Algebra aus; Hilbert hatte nat\"urlich
gute Gr\"unde f\"ur seine Wahl) geh\"ort wohl auch zu den Punkten, die
Emmy Noether wiederholt zum Anlass f\"ur Kritik an Hilberts Zahlbericht
nahm (vgl. \cite{LS}).

Damit beansprucht der Prager Satz jedenfalls eine historische Rolle:
mit der Abkehr vom Prager Satz begann die Absatzbewegung weg von
rechnerischen hin zu konzeptionellen Beweisen und die Durchdringung
des untersuchten Gegenstands mit abstrakten Mitteln, eine Bewegung,
die Dedekind initiierte, Emmy Noether zu ihrem Programm erhob und
die mit dem Schaffen Grothendiecks\footnote{W. Scharlau hat j\"ungst
eine lesenswerte Biographie Grothendiecks (\cite{Sch}) ver\"offentlicht.}
einen vorl\"aufigen H\"ohepunkt erreicht hat.

\section{Plus \c{c}a change \ldots}

Wie wir gesehen haben, gilt das Gau{\ss}sche Hebungslemma f\"ur
normierte Polynome \"uber Integrit\"atsringen $R$
genau dann, wenn $R$ ganz abgeschlossen ist. F\"ur das allgemeine
Gau{\ss}sche Hebungslemma ist das nicht der Fall: das Polynom
$$ 2X^2-5 = \textstyle 2(X - \frac{\sqrt{10}}2) (X + \frac{\sqrt{10}}2)
          = (\sqrt{2} \, X - \sqrt{5}\,) (\sqrt{2} \, X + \sqrt{5}\,) $$
faktorisiert \"uber $\Q(\sqrt{10}\,)$ (und \"uber $\Z[\sqrt{2},\sqrt{5}\,]$),
nicht aber \"uber dem ganz abgeschlossenen Ring $R = \Z[\sqrt{10}\,]$.
In der Tat: w\"are das Hebungslemma \"uber $R$ richtig, so
m\"ussten wir $2$ in $R$ so faktorisieren k{\"o}nnen (also
$2 = \alpha \beta$ mit $\alpha, \beta \in R$), dass die Faktoren
$\alpha(X - \frac{\sqrt{10}}2)$ und $\beta(X + \frac{\sqrt{10}}2)$
Koeffizienten in $R$ haben. Da aber $2$ in $\Z[\sqrt{10}\,]$
irreduzibel ist, ist $\alpha$ oder $\beta$ eine Einheit, und dann
w\"are $\frac{\sqrt{10}}2 \in R$: Widerspruch.

Man mag nun vermuten, dass der Grund f\"ur das Scheitern des
Gau{\ss}schen Hebungslemmas darin liegt, dass der Ring $\Z[\sqrt{10}\,]$
nicht faktoriell ist. Tats\"achlich ist dies gar nicht so weit von
der Wahrheit entfernt: die G\"ultigkeit des Hebungslemmas in $R[X]$
hat etwas damit zu tun, ob in $R$ der Euklidische Vierzahlensatz gilt.
Mittels der Gleichung $2 \cdot 5 = \sqrt{10} \sqrt{10}$ sieht man
schnell ein, dass das Element $2$ in $\Z[\sqrt{10}\,]$ zwar
irreduzibel, aber nicht prim (und damit auch nicht primal) ist.
Insbesondere ist $\Z[\sqrt{10}\,]$ nicht Rieszsch. Jetzt zeigen wir

\begin{prop}
Ist $R[X]$ Rieszsch, so gilt in $R[X]$ das Hebungslemma.
\end{prop}

\begin{bew}
Sei $F \in R[X]$ und $F = gh$ mit $g, h \in K[X]$; o.B.d.A. sei
$F \ne 0$. W\"ahle $b, c \in R \setminus \{0\}$
derart, dass $bg, ch \in R[X]$ ist, und setze $a = bc$. Dann ist
$aF = (bg)(ch)$ in $R[X]$, und da $R[X]$ Rieszsch ist,
existieren $m, n, r, s \in R[X]$ mit $a = mr$, $bg = nr$, $ch = ms$
und $F = ns$. Wegen $\grad a = 0$ m\"ussen $m, r \in R$ sein. Weiter
liegen $G = \frac{bg}r = n$ und $H = \frac{ch}m = s$ beide in $R[X]$,
und es ist $F = GH$ mit $G = \frac br g$, $H = \frac cm h$ und
$\frac br \cdot \frac cm = 1$. Also gilt das Hebungslemma.
\end{bew}

Als n\"achstes m{\"o}chten wir die Umkehrung beweisen: Gilt in
$R[X]$ das Gau{\ss}sche Hebungslemma, so ist $R[X]$ Rieszsch.
Dazu m\"ussen wir zeigen, dass alle Elemente von $R[X]$ primal
sind. Ein erster Schritt in diese Richtung w\"are der Nachweis,
dass alle Elemente von $R$ primal in $R[X]$ sind. Ein noch
kleinerer Schritt schlie{\ss}lich ist das folgende

\begin{lem}\label{LeGLR}
Gilt in $R[X]$ das Gau{\ss}sche Hebungslemma, so ist $R$ Rieszsch.
\end{lem}

\begin{bew}
Sei $ad = bc$ f\"ur $a, b, c, d \in R$. Im Falle von $abcd = 0$ ist
die Behauptung leicht nachzuweisen; wir d\"urfen daher $abcd \ne 0$
annehmen. Jetzt betrachten wir das Polynom
$$ F(X) = a \Big(X - \frac ba\Big) \Big(X - \frac ca\Big)
        = aX^2 - (b+c)X + d \in R[X]. $$
Da wir diese Faktorisierung nach $R[X]$ hochheben k{\"o}nnen,
gilt $F(X) = G(X)H(X)$ f\"ur $G, H \in R[X]$, wobei $G$ ein konstantes
Vielfaches von $X - \frac ba$ und $H$ eines von $X - \frac ca$ ist.
Sei daher $G(X) = \alpha X - \beta$ und $H(X) = \delta X - \gamma$ mit
$\alpha, \beta, \gamma, \delta \in R$. Da $G$ und $g(X) = X - \frac ba$
dieselbe Wurzel haben, ist $\frac ba = \frac \beta\alpha$, also
$\alpha b = \beta a$; analog finden wir $\gamma a = \delta c$.
Weiter liefert Koeffizientenvergleich $a = \alpha \delta$ und
$d = \beta \gamma$. Daraus folgt
$bc = ad = a \beta \gamma = \beta\delta c$, also $b = \beta \delta$.
\"Ahnlich ergibt sich
$bc = a \beta \gamma = b\alpha \gamma$ und damit $c = \alpha \gamma$.
Zusammenfassend haben wir jetzt $a = \alpha \delta$, $b = \beta\delta$,
$c = \alpha\gamma$ und $d = \beta\gamma$, und das war zu zeigen.
\end{bew}

\noindent{\bf Bemerkung.} Man beachte, dass man f\"ur den Beweis von
Lemma \ref{LeGLR} nur das Gau{\ss}sche Hebungslemma f\"ur quadratische
Polynome vorauszusetzen braucht.
\medskip

Wir wissen jetzt, dass in Rieszschen Ringen $R[X]$ das Gau{\ss}sche
Lemma gilt, und dass aus der G\"ultigkeit des Gau{\ss}schen Lemmas
in $R[X]$ folgt, dass $R$ Rieszsch ist. Jetzt w\"urden wir gerne
zeigen, dass mit $R$ auch $R[X]$ Rieszsch ist; das geht
allerdings nicht:

\begin{prop}\label{PRX}
Sind alle $a \in R$ auch in $R[X]$ primal, dann ist $R$ ganz
abgeschlossen. Insbesondere ist $R$ ganz abgeschlossen, wenn
$R[X]$ Rieszsch ist.
\end{prop}

\begin{bew}
Sei $t \in K$ ganz \"uber $R$, also Nullstelle eines normierten Polynoms
$f \in R[X]$. Dann ist $f(X) = (X-t)g(X)$ f\"ur ein normiertes
$g \in K[X]$. Jetzt gibt es $r, s \in R \setminus \{0\}$ mit
$r(X-t), sg(X) \in R[X]$. Das Element
$q = rs \in R$ ist nach Annahme primal in $R[X]$, und es gilt
$q \mid [r \cdot (X-t)][s \cdot g(X)] = rs f(X)$. Also ist $q = ab$ mit
$a \mid r \cdot (X-t)$ und $b \mid s \cdot g(X)$, wobei $a, b \in R$
sind, da sie die Konstante $q$ teilen. Da $X-t$ und $g(X)$
normiert sind, muss $a \mid r$ und $b \mid s$ gelten. Wegen $ab = rs$
folgt daraus, dass $\frac ra, \frac sb \in R^\times$ Einheiten sind.
Aus $\frac ra(X-t) \in R[X]$ folgt dann aber $t \in R$, folglich
ist $R$ ganz abgeschlossen.
\end{bew}

Um daher die Rieszeigenschaft von $R$ auf $R[X]$ \"ubertragen
zu k{\"o}nnen, ben{\"o}tigen wir die ganze Abgeschlossenheit von $R$.
Wie sich nun herausstellt, brauchen wir diese f\"ur den Schluss vom
Hebungslemma auf die Rieszeigenschaft von $R[X]$ aber gar nicht
vorauszusetzen, da sie aus dem Hebungslemma folgt:

\begin{prop}\label{PGA}
Sei $R$ ein Integrit\"atsring. Gilt in $R[X]$ das Gau{\ss}sche
Hebungslemma, dann ist $R$ ganz abgeschlossen.
\end{prop}

\begin{bew}
Sei $K$ der Quotientenk{\"o}rper von $R$ und $t \in K$ ganz \"uber $R$,
also Nullstelle eines Polynoms
$f(X) = X^n + a_{n-1}X^{n-1} + \ldots + a_0 \in R[X]$. Dann gibt
es in $K[X]$ eine Faktorisierung $f(X) = (X-t)g(X)$ mit
$X-t, g(X) \in K[X]$. Nach dem Gau{\ss}schen Hebungslemma existiert
ein $q \in K^\times$ derart, dass $q(X-t)$ und $\frac1q g(X)$ beide
in $R[X]$ liegen. Da sowohl $X-t$ wie auch $g(X)$ normiert sind,
also Leitkoeffizient $1$ besitzen, muss $q \in R$ und $\frac1q \in R$
sein. Also ist $q \in R^\times$ Einheit, und dies impliziert
$X-t, g(X) \in R[X]$. Also ist $t \in R$, und $R$ ganz abgeschlossen.
\end{bew}

Insbesondere gilt das Gau{\ss}sche Hebungslemma nicht \"uber
$R = \Z[\sqrt{-3}\,]$: die Faktorisierung $X^2 + X + 1 = (X-\rho)(X-\rho^2)$
mit $\rho = \frac{-1+\sqrt{-3}}2 \in K = \Q(\sqrt{-3}\,)$ l\"asst sich
offenbar nicht zu einer Faktorisierung in $R[X]$ hochheben.

\section{Cohns Permanenzsatz}\label{Ab12}

Nun zur\"uck zu unserem Problem: wir wollen nachweisen, dass aus
der G\"ultigkeit des Gau{\ss}schen Hebungslemmas in $R[X]$ folgt, dass
$R[X]$ Rieszsch ist: wir haben bereits gesehen, dass $R$ dann Rieszsch
und ganz abgeschlossen ist; daraus m{\"o}chten wir gern schlie{\ss}en,
dass auch $R[X]$ Rieszsch ist. Der erste Schritt in diese Richtung
ist die folgende

\begin{prop}\label{PprimRX}
Sei $R$ Rieszsch und ganz abgeschlossen. Dann sind alle $a \in R$
auch primal in $R[X]$.
\end{prop}

\begin{bew}
Sei $a \in R \setminus \{0\}$ und $a \mid b(X)c(X)$ f\"ur
Polynome $b, c \in R[X] \setminus \{0\}$. Mit Dedekinds Prager Satz
folgt daraus\footnote{F\"ur Elemente $a, b, c \in R$ ist die Relation
$a \mid bc$ in $R$ gleichbedeutend mit derjenigen in $R[X]$,
sodass wir zwischen $a \mid_R bc$ und $a \mid_{R[X]} bc$ nicht zu
unterscheiden brauchen.} $a \mid b_ic_j$ f\"ur alle $i, j$.
Da $R$ Rieszsch ist, ist die Teilbarkeitsgruppe $G(R)$ eine Rieszsche
Gruppe, und nach der Euklidischen Verfeinerungseigenschaft aus Teil I
(angewandt auf die von $0$ verschiedenen Koeffizienten; Koeffizienten
$0$ werden von jedem Element geteilt) gibt es eine Faktorisierung
$a = rs$ mit $r \mid b_i$ und $s \mid c_j$ f\"ur alle $0 \le i \le m$
und $0 \le j \le n$. Also ist $r \mid b(X)$ und $s \mid c(X)$; mit
anderen Worten: $a$ ist auch primal in $R[X]$.
\end{bew}

\noindent{\bf Bemerkung.} Hier gen\"ugt es nicht, nur die Primalit\"at
des einzelnen Elements $a$ in $R$ vorauszusetzen: die Euklidische
Verfeinerungseigenschaft $E_{m,n}$ benutzt die Interpolationseigenschaft
der ganzen Gruppe $G(R)$, nicht nur die Primalit\"at von $a$.
\medskip

Aus unseren obigen Ausf\"uhrungen sch\"alt sich nun folgender Satz
heraus:

\begin{thm}\label{SRX}
Sei $R$ ein Integrit\"atsring. Dann sind folgende Eigenschaften
\"aquivalent:
\begin{enumerate}
\item $R[X]$ ist Rieszsch.
\item In $R[X]$ gilt das Gau{\ss}sche Hebungslemma.
\item $R$ ist Rieszsch und ganz abgeschlossen.
\item $R[X]$ ist Rieszsch und ganz abgeschlossen.
\end{enumerate}
\end{thm}

Ganz abgeschlossene Rieszsche Ringe nennt man (nach Cohn) auch
Schreiersche Ringe.

Bisher haben wir bewiesen, dass (1) $\impl$ (2) $\impl$ (3) gilt, und
selbstverst\"andlich ist (4) $\impl$ (1) trivial. Es bleibt also noch
der ``Permanenzsatz'' (3) $\impl$ (4) f\"ur die Schreier-Eigenschaft
zu beweisen. Dass das Hebungslemma genau in Schreierschen Ringen gilt
(dies stammt von McAdam \& Rush \cite{MR}), ist also, da die andern
Implikationen ziemlich naheliegend sind, im Wesentlichen \"aquivalent zum
Cohnschen Permanenzsatz (3) $\impl$ (4) (vgl. dazu \cite[Thm. 9.12]{MM}).

Es gibt einige Eigenschaften, die sich von einem Ring $R$ auf den
Polynomring $R[X]$ vererben. So ist bekannt, dass mit $R$ auch immer
$R[X]$ faktoriell ist\footnote{Dieser Satz wird oft Satz von Gau{\ss}
genannt, obwohl Gau{\ss} ihn weder formuliert noch bewiesen hat.
Erstmals explizit formuliert und bewiesen hat den Satz wohl erst Hensel
\cite{Hen27}, obwohl nat\"urlich Spezialf\"alle wie die ZPE-Eigenschaft
der Polynomringe $\Z[X_1,\ldots,X_n]$ bereits Kronecker bekannt
waren (wom\"oglich tauchte der Satz auch schon in Algebravorlesungen
z.B. von Emil Artin und Emmy Noether auf; jedenfalls findet er sich
wenig sp\"ater in van der Waerdens Buch \cite{vdWA}, das auf
diesen Vorlesungen basierte). Die Bezeichnung ``Satz von Gau{\ss}''
r\"uhrt wohl von folgender Bemerkung van der Waerdens in \cite{vdWA} her:
``Der hier darzustellende Beweis geht auf Gau\ss\ zur\"uck''. In der
Tat f\"uhrt er den Beweis des Satzes mit Hilfe des Gau{\ss}schen Lemmas.}.
Auch die Eigenschaft, GGT-Ring (oder ganz abgeschlossen) zu sein,
\"ubertr\"agt sich von $R$ auf $R[X]$. Dass sich die Eigenschaft
Rieszsch nicht immer \"ubertr\"agt, folgt aus Prop. \ref{PRX} (und einem
Beispiel eines Rings, der Rieszsch, aber nicht ganz abgeschlossen ist).

W\"ahrend man f\"ur den ``Satz von Gau\ss'' einige einfache Beweise
kennt, scheint der Beweis von (3) $\impl$ (4) tiefer zu liegen; da
mir auch nur der Beweis Cohns\footnote{Der hier gegebene Beweis ist
anders strukturiert als der Cohnsche (dieser geht den Weg \"uber ein
dem Satz von Nagata (sh. z.B. \cite{LL}) analoges Ergebnis, welches
einen Zusammenhang zwischen der Schreiereigenschaft eines Rings $R$ und
seiner Lokalisierungen $R_S$ herstellt), benutzt aber im Wesentlichen
dieselben Hilfsmittel: Interpolations- bzw. Verfeinerungseigenschaften
in Rieszschen Ringen und Dedekinds Prager Satz.} \cite{Cohn} und
dessen Version von Anderson \& Zafrullah\footnote{Ihr Beweis benutzt
wie der unsere die Euklidische Verfeinerungseigenschaft.} \cite{AZ}
bekannt sind, sei den Lesern hiermit ans Herz gelegt, sich weitere
und vielleicht direktere Beweise f\"ur diesen Satz zu \"uberlegen.

Um zu zeigen, dass $R[X]$ unter geeigneten Voraussetzungen
Rieszsch ist, m\"ussen wir aus $a \mid bc$ f\"ur Polynome
$a, b, c \in R[X]$ die Existenz einer Faktorisierung
$a = rs$ mit $r, s \in R[X]$ schlie{\ss}en, f\"ur die
$r \mid b$ und $s \mid c$ gilt. In GGT-Ringen k{\"o}nnten wir
einfach $r = \ggT(a,b)$ setzen; in unserem Fall geht das
auch, wenn wir von $R[X]$ zum euklidischen Ring $K[X]$ \"ubergehen,
wo $K$ der Quotientenk{\"o}rper von $R$ ist. Damit erhalten wir eine
Faktorisierung $a = \rho\sigma$ mit $\rho, \sigma \in K[X]$
und $\rho \mid b$, $\sigma \mid c$ ebenfalls in $K[X]$.
Um daraus eine Faktorisierung in $R[X]$ zu gewinnen, liegt es
nahe, das Gau{\ss}sche Hebungslemma anzuwenden. Wir zeigen daher
zuerst die Richtung (3) $\impl$ (2) aus Satz \ref{SRX}:

\begin{prop}
Ist $R$ Rieszsch und ganz abgeschlossen, so gilt in $R[X]$
das Gau{\ss}sche Hebungslemma.
\end{prop}

\begin{bew}
Sei $f \in R[X]$ und $f = gh$ mit $g, h \in K[X] \setminus \{0\}$. Dann
existieren $b, c \in R \setminus \{0\}$ mit $bg, ch \in R[X]$, und mit
$a = bc$ ist $af = (bg)(ch)$, also $a \mid (bg)(ch)$ in $R[X]$. Nach
dem Prager Satz ist $a \mid (bg)_i(ch)_j$ f\"ur alle $i, j$, wo
$(bg)_i = bg_i$ und $(ch)_j = ch_j$ die Koeffizienten von $bg$ bzw.
$ch$ bezeichnen. Die Euklidische Verfeinerungseigenschaft in der
Rieszschen Gruppe $G(R)$ gibt uns eine Faktorisierung $a = rs$ mit
$r \mid bg_i$ und $s \mid ch_j$ f\"ur alle $i, j$. Damit sind
$G = \frac br g$ und $H = \frac cs h$ Polynome in $R[X]$, und es gilt
$GH = \frac{bc}{rs} gh = f$. Das war zu zeigen.
\end{bew}

Jetzt versch\"arfen wir die Aussage von Lemma \ref{LeGLR}
und zeigen (2) $\impl$ (1) aus Satz \ref{SRX}; da wir bereits
wissen, dass aus der G\"ultigkeit des Gau{\ss}schen Hebungslemmas
in $R[X]$ die ganze Abgeschlossenheit von $R$ und damit auch die
von $R[X]$ folgt, haben wir damit auch (2) $\impl$ (4) und somit den
ganzen Satz \ref{SRX} bewiesen.

\begin{prop}
Gilt in $R[X]$ das Gau{\ss}sche Hebungslemma, so ist $R[X]$ Rieszsch.
\end{prop}

\begin{bew}
Nach Lemma \ref{LeGLR} und Prop. \ref{PGA} ist $R$ Rieszsch und ganz
abgeschlossen; insbesondere gilt in $R[X]$ Dedekinds Prager Satz
(Satz \ref{GLDP}), und die Teilbarkeitsgruppe $G(R) = K^\times/R^\times$
ist eine Rieszsche Gruppe, besitzt also die Interpolationseigenschaft
$R_{m,n}$ f\"ur alle $m, n$ (Teil I). Wir schreiben wieder $x \le y$
f\"ur $x, y \in K^\times$, wenn $\frac yx \in R$ ist; insbesondere
bedeutet $1 \le x$ einfach $x \in R$.

Sei $a \mid bc$, also $ad = bc$ f\"ur Polynome $a, b, c, d \in R[X]$.
Ist eines dieser Polynome gleich $0$, so ist die Behauptung leicht
nachzupr\"ufen; seien also $a, b, c, d \ne 0$.
Da $K[X]$ ein GGT-Ring ist, gibt es Polynome $r, s, t, u \in K[X]$
mit $a = rs$, $d = tu$, sowie $b = rt$ und $c = su$.
Weil in $R[X]$ das Gau{\ss}sche Hebungslemma gilt, d\"urfen wir
$r, s, t, u$ so um konstante Faktoren ab\"andern, dass
$r, s, t, u \in R[X]$ wird und die Gleichungen $a = rs$ und $d = tu$
erhalten bleiben. Die beiden andern Gleichungen gelten dann nur noch
bis auf konstante Faktoren: es gibt ein $q \in K^\times$ mit $b = q\,rt$,
und wegen $rt \cdot su = ad = bc = q\,rt \cdot c$ ist dann
$c = q^{-1} su$.

Schreiben wir $q = \frac\beta\alpha$ mit $\alpha, \beta \in R$, so ist
$\alpha b = (\beta r) t$, also $\alpha \mid (\beta r)t$ in $R[X]$, und
Dedekinds Prager Satz liefert $\alpha \mid \beta r_it_j$, d.h.
$1 \le q r_i t_j$ f\"ur alle von $0$ verschiedenen Koeffizienten $r_i$
von $r$ und $t_j$ von $t$. Ebenso folgt aus $\beta c = (\alpha s)u$,
dass $1 \le q^{-1} s_ku_l$ f\"ur alle von $0$ verschiedenen
Koeffizienten $s_k$ von $s$ und $u_l$ von $u$ gilt.

Wir suchen nun ein $\kappa \in K^\times$ derart, dass
$$ r' = \kappa r, \quad s' = \kappa^{-1} s, \quad
   t' = q\kappa^{-1} t, \quad u' = q^{-1}\kappa u $$
in $R[X]$ liegen. Dann ist n\"amlich $a = rs = r's'$,
$d = tu = t'u'$, $b = q\,rt = r't'$ und $c = q^{-1}su = s'u'$,
und $a \mid bc$ ist in der gew\"unschten Form zerlegt. Die Bedingungen
$r', s', t', u' \in R[X]$ lauten
$$ r_i^{-1} \le \kappa, \quad q u_l^{-1} \le \kappa, \quad
   \kappa \le s_k, \quad \kappa \le q t_j $$
f\"ur alle $i, j, k, l$ (mit von $0$ verschiedenen Koeffizienten).
Nach der Rieszschen Interpolationseigenschaft existiert ein solches
$\kappa$, sofern jede der linken Schranken unterhalb jeder der
rechten liegt. Das ist aber der Fall:
\begin{itemize}
\item $r_i^{-1} \le s_k$ bedeutet $1 \le r_is_k$, und das ist
      wegen $r_is_k \in R$ richtig;
\item $r_i^{-1} \le q t_j$ bedeutet $1 \le q r_i t_j$, s.o.;
\item $q u_l^{-1} \le s_k$ bedeutet $1 \le q^{-1} s_k u_l$, s.o.;
\item $q u_l^{-1} \le q t_j$ bedeutet $1 \le t_j u_l$, und das ist
      wegen $t_ju_l \in R$ richtig.
\end{itemize}
Damit ist der Beweis erbracht.
\end{bew}

Es sei noch einmal darauf hingewiesen, dass ein konzeptioneller
Beweis des Cohnschen Permanenzsatzes sehr willkommen w\"are.

\section{Andere Formen des Gau\ss{}schen Lemmas}\label{Ab13}

In vielen Quellen wird das Gau{\ss}sche Lemma anders formuliert.
Ist $R$ ein Integrit\"atsring, so nennt man ein Polynom $f \in R[X]$
{\em primitiv}, wenn jeder gemeinsame Teiler der Koeffizienten von $f$
eine Einheit in $R$ ist. Ist $R$ ein GGT-Ring und
$f(X) = f_0+f_1X + \ldots + f_nX^n$, so nennt man
$\ct(f) = \ggT(f_0, f_1, \ldots, f_n)$ den Inhalt von $f$.
Wie Dedekind (im Falle $R = \Z$; sh. \cite{Ded}) bemerkt hat,
l\"asst sich das Gau{\ss}sche Hebungslemma dann wie folgt formulieren.

\begin{lem}\label{Lggt}
Ist $R$ ein faktorieller Ring (oder auch nur ein GGT-Ring), so ist
das Produkt primitiver Polynome in $R[X]$ wieder primitiv. Genauer
gilt: f\"ur $f, g \in R[X]$ ist $\ct(fg) = \ct(f) \ct(g)$ (bis auf
Einheiten).
\end{lem}

Der Inhalt $\ct(f)$ von Polynomen hat die Eigenschaft $\ct(ag) = a \ct(g)$
f\"ur $a \in R$ und $g \in R[X]$; insbesondere folgt aus $a \mid f$
f\"ur $f \in R[X]$ immer $a \mid \ct(f)$.

Das Hebungslemma f\"ur GGT-Ringe ergibt sich damit so: ist $F = gh$
f\"ur Polynome $F \in R[X]$ und $g, h \in K[X]$, so gibt es
$a, b \in R \setminus \{0\}$ mit $ag, bh \in R[X]$. Damit ist
$abF = (ag)(bh)$ eine Zerlegung in $R[X]$, und es gilt
$ab\ct(F) = \ct(ag) \ct(bh)$. Mit $r = \ggT(ab,\ct(ag))$
und $s = \frac{ab}r$ ist dann $ab = rs$ mit $r \mid ag$; da
$s = \frac{ab}r$ und $\frac{\ct(ag)}r$ teilerfremd sind und
$s \mid \frac{\ct(ag)}r \ct(bh)$ gilt, folgt auch $s \mid \ct(bh)$,
also $s \mid bh$.
Also sind $G = \frac ar g$ und $H = \frac bs h = \frac ra h$ Polynome
in $R[X]$ mit $F=GH$.

\medskip\noindent{\bf Definition.}
Wir nennen $R$ einen $P^4$-Ring (f\"ur $P^4 = PPPP$: Produkte primitiver
Polynome sind primitiv), wenn das Produkt zweier primitiver Polynome in
$R[X]$ wieder primitiv ist.
\medskip

Wir haben gesehen, dass das Gau{\ss}sche Hebungslemma f\"ur
normierte Polynome in $R[X]$ genau dann gilt, wenn $R$
(und damit $R[X]$) ganz abgeschlossen ist, bzw. wenn
Dedekinds Prager Satz gilt. Entsprechend gilt das
Gau{\ss}sche Hebungslemma f\"ur beliebige Polynome genau in
Schreierschen Ringen. Es ist auch nicht schwer zu zeigen,
dass $P^4$-Ringe zwischen Rieszsch und AP liegen (sh. Abb.
\ref{FR}; tats\"achlich folgt AP bereits aus der G\"ultigkeit
der $P^4$-Eigenschaft f\"ur lineare Polynome). Das
Gau{\ss}sche Lemma in der Form von Lemma \ref{Lggt} gilt
schlie{\ss}lich genau dann, wenn $R$ ein GGT-Ring ist.

Nun kann man in beliebigen Integrit\"atsringen das von den
Koeffizienten von $f$ erzeugte Ideal $\fc(f)$ betrachten und
sich fragen, unter welchen Voraussetzungen die Gleichung
$\fc(fg) = \fc(f)\fc(g)$ gilt. Die Antwort darauf gibt
folgender

\begin{thm}
F\"ur einen Integrit\"atsring $R$ sind folgende Aussagen \"aquivalent:
\begin{enumerate}
\item Jedes endlich erzeugte Ideal ist projektiv.
\item Jedes endlich erzeugte Ideal $\ne (0)$ ist invertierbar.
\item F\"ur jedes Primideal $P \ne (0)$ ist die Lokalisierung $R_P$ ein
      Bewertungsring.
\item Jeder Ring $S$ mit $R \subseteq S \subseteq K$, wo $K$ der
      Quotientenk{\"o}rper von $R$ ist, ist ganz abgeschlossen.
\item Es ist $\fc(fg) = \fc(f)\fc(g)$ f\"ur alle Polynome $f, g \in R[X]$.
\end{enumerate}
\end{thm}

\begin{bew}
Sh. \cite[Kap. IV]{Gil}.
\end{bew}

Diese Version des Gau{\ss}schen Lemmas gilt also genau dann,
wenn der Koeffizientenring Pr\"ufersch ist. Pr\"ufer selbst
(sh. \cite[S. 23]{Pru}) hat den Satz $\fc(fg) = \fc(f)\fc(g)$
den ``Satz von Gau\ss{}-Kronecker-Dedekind'' genannt.

\section*{Zu guter Letzt \ldots}
Die Motivation zum Schreiben dieses Artikels r\"uhrt im Wesentlichen
von meiner Lekt\"ure von \cite{MM} (darin geht es um das Gau{\ss}sche
Hebungslemma in Schreierschen Ringen) und \cite{PR} (wo die L\"ucke
in Euklids Beweis seiner Prop. VII.19 thematisiert wird) und dem Wunsch
her, das dort Gesagte ein wenig anders aufzubereiten. Weiter m{\"o}chte
ich Bill Dubuque f\"ur diverse Beitr\"age in sci.math danken, die mir vor
allem zu Anfang sehr hilfreich waren. Endlich verdanke ich Heinz
L\"uneburg und Klaus Volkert sowie den beiden anonymen Referenten
noch zahlreiche n\"utzliche und hilfreiche Bemerkungen.

\newpage

\begin{otherlanguage*}{british}
  
\renewcommand{\proofname}{Proof}
\renewcommand{\theHsection}{en.\arabic{section}}

\begin{center}
{\Large Part II. Gauss's Lemmas and Riesz Domains}
\end{center}
\medskip

\medskip\noindent
\begin{quote}
  \begin{small}
    \textsc{Abstract.}
    In the second part (cf.\ \cite{LSb}) of our article on Euclid's
    four number theorem we start from the notion of incommensurability
    and take a closer look at irrationality proofs for algebraic
    numbers.  Along the way we encounter Gauss's Lemma, integral
    closure and Dedekind's Prague theorem, and finally we arrive back
    at the four number theorem.
  \end{small}  
\end{quote}

\setcounter{section}{7}

\section{The Ancient Greeks once more: incommensurability}

As is well known, the Pythagoreans are also credited with the
discovery of incommensurability, i.e., with the observation that,
e.g., the diagonal $d$ and the side $s$ of a square do not stand in a
rational ratio, so that one cannot write $d:s = p:q$ with natural
numbers $p, q$. There is a multitude of proofs of this claim, two of
which go back to the Greeks:
\begin{itemize}
\item One proof deduces from $\sqrt{2} = \frac pq$ that $p$ and $q$
      must both be even, which contradicts the reducibility of
      fractions (or proportions). This proof goes back to Aristotle
      and is contained in some (older) editions of the Elements as
      Proposition X.117; today, however, it is regarded as a later
      interpolation.
\item The other proof (cf.\ Hasse \cite{Has} and Apostol \cite{Apo})
      is geometric in nature and is based on ``anthyphairesis''
      (reciprocal subtraction), the father of the Euclidean
      algorithm.
\end{itemize}
That the Pythagoreans already discovered the existence of
incommensurable segments is considered certain; whether it was
Hippasus who made (or made public) this discovery is uncertain, as is
the claim that he was led to it by studying the pentagram (see von
Fritz \cite{vFr}). Finally, the occasionally encountered claim that
this caused a ``foundational crisis'' of Greek mathematics can no
longer be supported by any historical facts (which does not mean that
it must be false).

In the course of the extension of the concept of number, the
incommensurability of segments has turned from a geometric into an
arithmetic problem: for us, the ratio of two segments is simply a
number, and the incommensurability of side and diagonal of a square
is the irrationality of the square root of $2$.

If one knows that the ring $\Z$ is factorial, the proof of the
irrationality of $\sqrt{2}$ becomes trivial: from $2q^2 = p^2$ one
reads off that $2$ divides the left hand side to an odd power and the
right hand side to an even power: contradiction.

The same proof shows more generally\footnote{Without unique
factorization this goes wrong: in the monoid $M$ of all natural
numbers $\equiv 0, 1, 2 \bmod 4$, the number $9$ is the square of the
fraction $\frac62$ in the quotient monoid, but not the square of an
element of $M$. Similarly, $2i$ is not a square in
$\Z[2i] = \Z \oplus 2i\Z$, but it is a square in the quotient field
$\Q(i)$ because $2i = (1+i)^2 = (\frac{2+2i}{2})^2$.}:

\begin{thm}\label{eSm}
For $m \in \Z$, $\sqrt{m}$ is rational if and only if $m$ is the
square of an integer.
\end{thm}
This, too, is a theorem that was essentially known to the Greeks; it
cannot, however, be found in this explicit form in Euclid's Elements.

Instead of proving the irrationality of $\sqrt{2}$ with the help of
unique prime factorization, we can also obtain it using Euclidean
division with remainder (which, in view of the fact that Euclidean
rings are factorial, should of course come as no surprise): let
$m \in \N$ be a non-square; then there is an $n \in \N_0$ with
$n^2 < m < (n+1)^2$, i.e., with $0 < \sqrt{m}-n < 1$. If $\sqrt{m}$
is rational, we can write $\sqrt{m} = \frac pq$ with minimal
$q \in \N$. Then $r := q(\sqrt{m}-n) = p - nq$ is an integer as well,
and $0 < \sqrt{m}-n < 1$ implies $0 < r < q$. Now
$r\sqrt{m} = q(\sqrt{m}-n)\sqrt{m} = mq-nq\sqrt{m} = mq - np$ is an
integer; since $\sqrt{m} = \frac{r\sqrt{m}}{r}$, this contradicts the
minimality of $q$. For generalizations of this kind of
argument\footnote{This proof essentially goes back to Dedekind
\cite{DedSIZ}; for a geometric version of the proof due to Haj\'os
see \cite{Sur2}. Variants of this proof can be found, among others,
in \cite{Blo,Fine,Sag,Vau,Wat}.} as well as other irrationality
proofs see, e.g., the exercises in \cite[Ch. 1]{LL}.

Here is a proof based on the fact that $\Z$ is a Bezout domain; it
shows more generally that a rational $n$-th root of an integer $m$ is
already an integer: let $\sqrt[n]{m} = \frac pq$ with coprime
integers $p$ and $q$, i.e., $(p,q) = (1)$. Then $p^n = mq^n$ gives
$$ (1) = (p,q)^n \subseteq (p^n,q) = (mq^n,q) \subseteq (q), $$
hence $(q) = (1)$, so $q = \pm 1$ and $\sqrt[n]{m} = \pm p$ is an
integer.

As pretty as such proofs may be, in a certain sense they obscure the
actual core of the matter; we want to uncover it in the next section.

\section{Integral closure}

Students usually encounter the notion of integral closure in
introductory algebra courses: if $R$ is an integral domain with
quotient field $K$, an element $t \in K$ is called integral over $R$
if $t$ is a root of a monic polynomial $f \in R[X]$. If every
$t \in K$ that is integral over $R$ already lies in $R$, then $R$ is
called integrally closed (in its quotient field).

If $t_1, t_2$ (in some extension field of $K$) are integral over
$R$, then so are $t_1 + t_2$ and $t_1t_2$. The classical proof using
resultants has been replaced over time by proofs via finitely
generated modules; dressed in the language of tensors, however, the
classical proof almost matches the module proof in elegance: as a
root of a monic polynomial of degree $m_j$ with coefficients in $R$,
$t_j$ is also an eigenvalue of an $m_j \times m_j$ matrix $T_j$ (the
companion matrix) with entries in $R$. One now easily shows that
$t_1t_2$ is an eigenvalue of the Kronecker product $T_1 \otimes T_2$
and hence (as a root of its characteristic polynomial) integral over
$R$. Similarly, $t_1 + t_2$ is an eigenvalue of
$T_1 \otimes I_{m_2} + I_{m_1} \otimes T_2$, where $I_k$ denotes the
identity matrix of dimension $k$ (see \cite{KMS}).

The theorem that factorial domains are integrally closed belongs to
the basic knowledge of algebra students; the usual proof even gives
the stronger statement that GCD domains are integrally closed. But
even in courses on algebraic number theory it is rarely made clear
how fundamental the notion of integral closure has been for the
development of this discipline, and what significance it has gained
for modern algebra since Dedekind and Emmy Noether
distilled\footnote{Kummer did not know the notion of an algebraic
integer, and the significance of integral closure did not become
clear to him. His construction of ideal numbers in cyclotomic fields
worked because the rings Kummer considered were integrally closed;
his development of the theory of ideal numbers, however, has a gap
that was only closed by Dedekind -- using integral closure (cf.\
\cite{LJ}).

Dirichlet and Eisenstein did not see integral closure either:
Dirichlet proved his unit theorem in rings of the form $\Z[\alpha]$,
where $\alpha$ is an algebraic integer, and it follows from a remark
of Eisenstein that the algebraic integers form a ring. Heine
\cite{Heine} defined algebraic integers as those numbers that can be
obtained from the integers by ring operations and repeated extraction
of $n$-th roots.

Dedekind showed that in the ring of integers of a number field
$K \subseteq \C$ every ideal $\ne (0)$ is uniquely a product of prime
ideals; finally, in 1926 Emmy Noether succeeded in showing that this
theorem holds in all integral domains satisfying certain axioms --
one of which is integral closure.} it from the construction of
Dedekind's ideal theory. Indeed, a look at the following figure shows
that most of the classes of rings\footnote{Notions not yet defined,
such as Schreier or $P^4$, will be explained in Sections \ref{eAb12}
and \ref{eAb13}.} occurring in introductory algebra courses are
integrally closed:

\begin{figure}[h!]
$$ \begin{array}{cccccccccccc}
\text{Euclid} &         &       &    & & & \text{Dedekind} & & \\
    \cap      &         &       &    & & & \cap & & \\
 \text{PID}    & \subset & \text{Bezout} & & \subset & & \text{Pr\"ufer} & & \\
    \cap      &         &   \cap     &         & & & \cap & & \\
\text{UFD} & \subset & \text{GCD} & \subset & \text{Schreier}
              & \subset & \text{int. closed} &  \\
    \cap      &         &      &       &  \cap   &  &  & & \\
\text{atomic} &         &      &   & \text{Riesz} & \subset & P^4 &
    \subseteq & \text{AP} \\
\end{array}$$
\caption{Classes of rings}\label{eFR}
\end{figure}

The integral closure of $\Z$ is also the real reason for the validity
of Theorem \ref{eSm}. Indeed, if $n \in \Z$ is the square of a
rational number, say $\sqrt n = \frac pq$, then $\frac pq$ is a root
of the quadratic polynomial $X^2 - n \in \Z[X]$. Since $\Z$ is
integrally closed, we must have $\frac pq \in \Z$.

It is therefore not surprising at all that the irrationality of
$\sqrt{2}$ can be proved using the factoriality of $\Z$. Moreover,
one can try to give such irrationality proofs using any of the
theorems from which the integral closure of $\Z$ follows as well. An
example is the following proof:

\begin{align*}
 \sqrt{m} \text{ is rational}
     & \iff X^2 - m \text{ is reducible in } \Q[X] \\
     & \iff X^2 - m \text{ is reducible in } \Z[X] \\
     & \iff X^2 - m = (X-n)(X+n) \text{ for some } n \in \Z \\
     & \iff m = n^2 \text{ for some } n \in \Z.
\end{align*}
The passage from $\Q[X]$ to $\Z[X]$ is achieved by Gauss's lifting
lemma\footnote{The name is meant to indicate that factorizations of
polynomials over $\Q$ can be lifted to factorizations over $\Z$.} for
monic polynomials:

\begin{lem}\label{eGHL1}
If $g, h \in \Q[X]$ are monic polynomials with $gh = f \in \Z[X]$,
then $g, h \in \Z[X]$.
\end{lem}

Gauss's lifting lemma appears in this form in the Disquisitiones
Arithmeticae (namely in \cite[Art. 42]{DA}), where it is used to
prove the irreducibility of the cyclotomic polynomials
$F_p(X) = \frac{X^p-1}{X-1}$ for primes $p$. In order to distinguish
it from other results that are also known as ``Gauss's
Lemma''\footnote{There are at least three such results: the lifting
lemma above, and one each from the theory of quadratic residues and
the theory of quadratic forms. All three can be found in the
Disquisitiones Arithmeticae.

Moreover, Gauss's lifting lemma in $\Z$ is often formulated
differently, and while these versions are equivalent to the lifting
lemma formulated above over $\Z$, they are not over general integral
domains.}, we will call Lemma \ref{eGHL1} Gauss's lifting lemma (for
monic polynomials).

Lemma \ref{eGHL1} allows us to show much more than just the
irrationality of certain square roots: if $\alpha$ is a root of a
monic polynomial with integral coefficients (i.e., an algebraic
integer), then $\alpha$ is either an integer or irrational. In fact,
if $\alpha$ is a rational root of the monic polynomial
$f \in \Z[X]$, then $f(X) = (X-\alpha)g(X)$ with $g \in \Q[X]$ monic.
Gauss's lifting lemma for monic polynomials then shows
$X - \alpha \in \Z[X]$, hence $\alpha \in \Z$.

Gauss's lifting lemma for monic polynomials is a special case of

\begin{thm}[Gauss's lifting lemma]\label{eLGHL}
If $F \in \Z[X]$ is a polynomial and if there are polynomials
$g, h \in \Q[X]$ with $F = gh$, then there is a $c \in \Q^\times$
such that $G = cg, H = \frac1c h \in \Z[X]$. Thus every factorization
$F = gh$ in $\Q[X]$ can be ``lifted'' to a factorization $F = GH$ in
$\Z[X]$.
\end{thm}

This theorem is still an important tool in today's algebra courses
for proving the irreducibility of polynomials and is used, for
example, in the proof of Eisenstein's criterion: with its help one
shows that a primitive polynomial of degree $\ge 1$ in $\Z[X]$ is
irreducible in $\Q[X]$ if and only if it is irreducible in $\Z[X]$.

If $R$ is an integral domain with quotient field $K$, we say that
Gauss's lifting lemma holds in $R[X]$ if every factorization of an
$F \in R[X]$ can be lifted from $K[X]$ to $R[X]$, i.e., if for every
factorization $F = gh$ with $g, h \in K[X]$ there is a
$c \in K^\times$ with $cg, \frac1c h \in R[X]$. We now ask in which
rings the above versions of Gauss's lifting lemma hold. For monic
polynomials this question is relatively easy to answer:

\begin{thm}\label{eGLDP}
Let $R$ be an integral domain with quotient field $K$. Then the
following statements are equivalent:
\begin{enumerate}
\item $R$ is integrally closed.
\item Gauss's lifting lemma for monic polynomials holds in $R[X]$: \\
      if $g, h \in K[X]$ are monic polynomials with $gh = f \in R[X]$,
      then $g, h \in R[X]$.
\item Dedekind's Prague theorem holds in $R[X]$: if \\
      \begin{equation}\label{eEDPS}
       b(X) = \sum_{i=0}^m b_i X^i \quad \text{and } \quad
       c(X) = \sum_{j=0}^n c_jX^j
      \end{equation}
      are polynomials in $R[X]$ and $a \mid b(X)c(X)$ for some
      $a \in R$, then $a \mid b_ic_j$ for all $0 \le i \le m$,
      $0 \le j \le n$.
\end{enumerate}
\end{thm}

Before we give the proof of Theorem \ref{eGLDP} and then study the
validity of the general lifting lemma, we make a few remarks on the
significance of Dedekind's Prague theorem.

\section{Dedekind's Prague theorem}

Historically, Dedekind's Prague theorem\footnote{Dedekind published
this theorem in the journal of the Prague Mathematical Society, and
today the name ``Prague theorem'' has become established for it.}
played an important role in the construction of Dedekind's ideal
theory. Dedekind felt that arguing with polynomials in the
construction of his ideal theory was inappropriate\footnote{Cf.\
Dedekind's comments on Hurwitz's approach to ideal theory in
\cite{DedH}, e.g.\ on p.~53, where he says that the use of functions
of variables always appeared to him as an aid foreign to the subject,
and on p.~55, where he says that the proof mentioned above could not
yet satisfy him completely, since in his view the purity of the
theory is clouded by the intrusion of functions of variables.} and
looked for other ways. Hurwitz \cite{Hurw} then published an approach
to ideal theory that was based on Dedekind's Prague theorem (which
Hurwitz did not know at the time), whereupon Dedekind declared that
he knew this route, but had preferred his own approach for
methodological reasons. Hilbert's Zahlbericht \cite{HZb} and
Landau's\footnote{Landau \cite[p. 145]{Lan} calls Dedekind's Prague
theorem a ``theorem of Kronecker''; in Cohn \cite{CAN} it is called
``Hurwitz's Lemma''. Walter Strobl points out in his thesis
\cite[pp. 157--161]{Str} that Bachmann discusses the Prague theorem
in detail in vol.~III and vol.~V of his five-volume treatise on
number theory. See also the introduction in Edwards \cite{EdDiv}.}
textbook \cite{Lan} on algebraic numbers both use Hurwitz's
approach\footnote{Dedekind's theorem serves to prove that for every
ideal $\fa \ne (0)$ in a maximal order of a number field there is an
ideal $\fb \ne (0)$ such that $\fa\fb$ is a principal ideal.},
whereas today one rather works with localizations (which surely would
have met with Dedekind's approval).

Dedekind \cite{Ded} gave two proofs of his Prague theorem; the first
can be found in Landau \cite{Lan} (Hilbert does not give a proof in
his Zahlbericht \cite{HZb}, but refers to papers by Dedekind
\cite{Ded}, Hurwitz \cite{Hurw}, Kronecker \cite{Kron} and Mertens
\cite{Mert}), and the second proof leads to the Dedekind--Mertens
Lemma, whose history and further development, among others by Artin
and Rees, would deserve an article of its own (cf.\
\cite{Coqu,CDLQ}).

Here we want to give a third proof, which is based on the observation
that Dedekind's Prague theorem holds exactly in integrally closed
domains\footnote{That Emmy Noether knew this Theorem \ref{eGLDP}
becomes clear from her remarks on Dedekind's article \cite{DedH} in
his collected works.}; thus, via integral closure, the Prague theorem
is equivalent to Gauss's lifting lemma for monic polynomials. The
crucial step (1) $\impl$ (3) is based on a classical lemma that
essentially goes back to Kronecker:

\begin{lem}[Kronecker]\label{eLKron}
Let $R$ be an integral domain with quotient field $K$, and let
$f(X) = \lambda \prod_{k=1}^N (X - \theta_k) \in R[X]$ with
$\lambda \ne 0$ and $\theta_k$ in some extension field $L$ of $K$.
Then for every subset $I \subseteq \{1, \ldots, N\}$ the element
$\lambda \theta_I$, where $\theta_I = \prod_{k \in I} \theta_k$, is
integral over $R$.
\end{lem}

\begin{bew}
Let $e_1, \ldots, e_N$ be the elementary symmetric functions of the
$\theta_k$; since $f = \lambda X^N - \lambda e_1 X^{N-1} \pm \ldots$,
we have $\lambda e_k \in R$ for all $k$, and of course
$\lambda \in R$. For fixed $r = |I|$ consider the polynomial
$$ P(Y) = \prod_{|J| = r} (Y - \lambda\theta_J), $$
where $J$ runs through all $r$-element subsets of
$\{1, \ldots, N\}$. The coefficient of $Y^{M-s}$ (with
$M = \binom Nr$) is $(-1)^s\lambda^s E_s$, where $E_s$ is the $s$-th
elementary symmetric function of the $\theta_J$. Now $E_s$ is a
symmetric polynomial in $\theta_1, \ldots, \theta_N$ with integer
coefficients in which each variable occurs at most to the $s$-th
power (each $\theta_J$ contains each $\theta_k$ at most once). By the
fundamental theorem on symmetric functions (in the sharper form
provided by its usual proof), $E_s$ is therefore a polynomial of
total degree $\le s$ in $e_1, \ldots, e_N$ with integer coefficients.
For every monomial $e_{k_1} \cdots e_{k_v}$ with $v \le s$, however,
$\lambda^s e_{k_1} \cdots e_{k_v} = \lambda^{s-v} (\lambda e_{k_1})
\cdots (\lambda e_{k_v}) \in R$. Thus $P$ is monic with coefficients
in $R$, and $\lambda\theta_I$, being a root of $P$, is integral over
$R$.
\end{bew}

\begin{bew}[of Theorem \ref{eGLDP}]
(1) $\impl$ (2): Let $L$ be a splitting field of $f$ over $K$; there,
$g$ and $h$ also split into linear factors, say
$g(X) = \prod_i (X - \alpha_i)$ and $h(X) = \prod_j (X - \beta_j)$.
Since $f(\alpha_i) = f(\beta_j) = 0$, the $\alpha_i$ and $\beta_j$ are
integral over $R$. Hence the elementary symmetric functions in the
$\alpha_i$ and in the $\beta_j$, i.e., the coefficients of $g$ and
$h$, are integral over $R$ as well. Since they lie in $K$ and $R$ is
integrally closed, they already lie in $R$, and we get
$g, h \in R[X]$ (this proof is taken from \cite{CS}).

(2) $\impl$ (1): Let $t \in K$ be integral over $R$, i.e., a root of
a monic polynomial $f \in R[X]$. Then $f(X) = (X-t)g(X)$ with a monic
$g \in K[X]$, and by (2) we have $X - t \in R[X]$, i.e., $t \in R$.

(1) $\impl$ (3): Let $a \mid b(X)c(X)$. If $a = 0$ or $b = 0$ or
$c = 0$, there is nothing to prove (if $a = 0$, then $bc = 0$). So
let $a, b, c \ne 0$, $b_m \ne 0 \ne c_n$, and
$f = \frac{bc}a \in R[X]$ with leading coefficient
$\lambda = \frac{b_mc_n}a$. In a splitting field $L$ of $bc$ over $K$
we have $b = b_m\prod_{k=1}^m (X-\beta_k)$,
$c = c_n \prod_{l=1}^n (X-\gamma_l)$, hence
$f = \lambda \prod_{k=1}^{m+n} (X-\theta_k)$ with
$(\theta_1, \ldots, \theta_{m+n}) = (\beta_1, \ldots, \beta_m,
\gamma_1, \ldots, \gamma_n)$. Now
$$ \frac{b_ic_j}a = \lambda \cdot \frac{b_i}{b_m} \cdot \frac{c_j}{c_n}
   = \pm \lambda \, e_{m-i}(\beta_1, \ldots, \beta_m) \,
           e_{n-j}(\gamma_1, \ldots, \gamma_n), $$
and expanding the right hand side gives a sum of terms of the form
$\lambda \theta_I$ with $I \subseteq \{1, \ldots, m+n\}$. By
Lemma \ref{eLKron}, $\frac{b_ic_j}a \in K$ is therefore integral over
$R$, hence lies in $R$, i.e., $a \mid b_ic_j$.

(3) $\impl$ (1): Let $t \in K$ be integral over $R$, i.e., a root of
a monic polynomial $f \in R[X]$. Then $f(X) = (X-t)g(X)$ with monic
$g \in K[X]$; let $a, b \in R \setminus \{0\}$ with $at \in R$ and
$bg \in R[X]$. Then $[a(X-t)][bg(X)] = abf(X)$ is divisible by $ab$,
and Dedekind's Prague theorem shows that $ab \mid atb$ in $R$ (up to
sign, the right hand side is the product of the constant term of
$a(X-t)$ and the leading coefficient $b$ of $bg(X)$); this implies
$t = \frac{atb}{ab} \in R$.
\end{bew}

\noindent{\bf Remark.} A direct argument from (2) to (3) of the kind
``$f = \frac{b}{a} \cdot c$ is a factorization in $K[X]$, hence it
can be lifted'' is not possible, since $b$ and $c$ are in general not
monic: for this one would need the general lifting lemma, and as we
shall see, this holds only in Schreier domains.\footnote{Paul Pollack
  pointed out the error in the original proof in 2016. The correction
  in this translation is due to Claude (anthropic).}
\medskip

Dedekind's Prague theorem plays practically no role in today's
algebra and number theory; as Dedekind himself remarked, this is only
consistent: it is equivalent to the integral closure of the base
ring, and instead of giving computational proofs with the Prague
theorem, one should, according to Dedekind, prefer conceptual proofs
using abstract principles such as integral closure. History has
proved him right, as Emmy Noether already wrote in her comments on
Dedekind's works, even though Hilbert in particular still preferred
Hurwitz's approach to ideal theory in his Zahlbericht. This ``wrong
choice'' (from the point of view of abstract algebra; Hilbert of
course had good reasons for his choice) presumably also belongs to
the points that Emmy Noether repeatedly took as an occasion to
criticize Hilbert's Zahlbericht (cf.\ \cite{LS}).

In any case, the Prague theorem can thus claim a historical role:
the turning away from the Prague theorem marked the beginning of the
movement away from computational and towards conceptual proofs, and
of the penetration of the objects under study by abstract means, a
movement which Dedekind initiated, which Emmy Noether made her
program, and which reached a provisional climax with the work of
Grothendieck\footnote{W. Scharlau has recently published a
biography of Grothendieck \cite{Sch} that is well worth reading.}.

\section{Plus \c{c}a change \ldots}

As we have seen, Gauss's lifting lemma for monic polynomials holds
over an integral domain $R$ if and only if $R$ is integrally closed.
This is not the case for the general lifting lemma: the polynomial
$$ 2X^2-5 = \textstyle 2(X - \frac{\sqrt{10}}2) (X + \frac{\sqrt{10}}2)
          = (\sqrt{2} \, X - \sqrt{5}\,) (\sqrt{2} \, X + \sqrt{5}\,) $$
factors over $\Q(\sqrt{10}\,)$ (and over $\Z[\sqrt{2},\sqrt{5}\,]$),
but not over the integrally closed ring $R = \Z[\sqrt{10}\,]$.
Indeed, if the lifting lemma held over $R$, we would have to be able
to factor $2$ in $R$ (i.e., $2 = \alpha \beta$ with
$\alpha, \beta \in R$) in such a way that the factors
$\alpha(X - \frac{\sqrt{10}}2)$ and $\beta(X + \frac{\sqrt{10}}2)$
have coefficients in $R$. But since $2$ is irreducible in
$\Z[\sqrt{10}\,]$, $\alpha$ or $\beta$ is a unit, and then
$\frac{\sqrt{10}}2 \in R$: contradiction.

One might now suspect that the reason for the failure of Gauss's
lifting lemma is that the ring $\Z[\sqrt{10}\,]$ is not factorial. In
fact, this is not that far from the truth: the validity of the
lifting lemma in $R[X]$ has something to do with whether Euclid's
four number theorem holds in $R$. Using the equation
$2 \cdot 5 = \sqrt{10} \sqrt{10}$ one quickly sees that the element
$2$ in $\Z[\sqrt{10}\,]$ is irreducible, but not prime (and hence not
primal either). In particular, $\Z[\sqrt{10}\,]$ is not a Riesz
domain. We now show

\begin{prop}
If $R[X]$ is a Riesz domain, then the lifting lemma holds in $R[X]$.
\end{prop}

\begin{bew}
Let $F \in R[X]$ and $F = gh$ with $g, h \in K[X]$; without loss of
generality $F \ne 0$. Choose $b, c \in R \setminus \{0\}$ such that
$bg, ch \in R[X]$, and put $a = bc$. Then $aF = (bg)(ch)$ in $R[X]$,
and since $R[X]$ is a Riesz domain, there exist
$m, n, r, s \in R[X]$ with $a = mr$, $bg = nr$, $ch = ms$ and
$F = ns$. Since $\deg a = 0$, we must have $m, r \in R$. Moreover,
$G = \frac{bg}r = n$ and $H = \frac{ch}m = s$ both lie in $R[X]$, and
$F = GH$ with $G = \frac br g$, $H = \frac cm h$ and
$\frac br \cdot \frac cm = 1$. Thus the lifting lemma holds.
\end{bew}

Next we would like to prove the converse: if Gauss's lifting lemma
holds in $R[X]$, then $R[X]$ is a Riesz domain. For this we have to
show that all elements of $R[X]$ are primal. A first step in this
direction would be to show that all elements of $R$ are primal in
$R[X]$. An even smaller step is the following

\begin{lem}\label{eLeGLR}
If Gauss's lifting lemma holds in $R[X]$, then $R$ is a Riesz domain.
\end{lem}

\begin{bew}
Let $ad = bc$ for $a, b, c, d \in R$. If $abcd = 0$, the claim is
easily verified; we may therefore assume $abcd \ne 0$. Now consider
the polynomial
$$ F(X) = a \Big(X - \frac ba\Big) \Big(X - \frac ca\Big)
        = aX^2 - (b+c)X + d \in R[X]. $$
Since we can lift this factorization to $R[X]$, we have
$F(X) = G(X)H(X)$ with $G, H \in R[X]$, where $G$ is a constant
multiple of $X - \frac ba$ and $H$ one of $X - \frac ca$. So write
$G(X) = \alpha X - \beta$ and $H(X) = \delta X - \gamma$ with
$\alpha, \beta, \gamma, \delta \in R$. Since $G$ and
$g(X) = X - \frac ba$ have the same root, we have
$\frac ba = \frac \beta\alpha$, i.e., $\alpha b = \beta a$; similarly
we find $\gamma a = \delta c$. Moreover, comparing coefficients gives
$a = \alpha \delta$ and $d = \beta \gamma$. This implies
$bc = ad = a \beta \gamma = \beta\delta c$, hence
$b = \beta \delta$. Similarly we get
$bc = a \beta \gamma = b\alpha \gamma$ and thus $c = \alpha \gamma$.
Summarizing, we now have $a = \alpha \delta$, $b = \beta\delta$,
$c = \alpha\gamma$ and $d = \beta\gamma$, which is what we wanted to
show.
\end{bew}

\noindent{\bf Remark.} Note that the proof of Lemma \ref{eLeGLR} only
requires Gauss's lifting lemma for quadratic polynomials.
\medskip

We now know that Gauss's Lemma holds in $R[X]$ if $R[X]$ is a Riesz
domain, and that the validity of Gauss's Lemma in $R[X]$ implies that
$R$ is a Riesz domain. We would now like to show that $R[X]$ is a
Riesz domain whenever $R$ is; this, however, is not possible:

\begin{prop}\label{ePRX}
If all $a \in R$ are also primal in $R[X]$, then $R$ is integrally
closed. In particular, $R$ is integrally closed if $R[X]$ is a Riesz
domain.
\end{prop}

\begin{bew}
Let $t \in K$ be integral over $R$, i.e., a root of a monic
polynomial $f \in R[X]$. Then $f(X) = (X-t)g(X)$ for some monic
$g \in K[X]$. Now there exist $r, s \in R \setminus \{0\}$ with
$r(X-t), sg(X) \in R[X]$. By assumption, the element $q = rs \in R$
is primal in $R[X]$, and $q \mid [r \cdot (X-t)][s \cdot g(X)] =
rsf(X)$. Thus $q = ab$ with $a \mid r \cdot (X-t)$ and
$b \mid s \cdot g(X)$, where $a, b \in R$ since they divide the
constant $q$. Since $X-t$ and $g(X)$ are monic, we must have
$a \mid r$ and $b \mid s$. Since $ab = rs$, this implies that
$\frac ra, \frac sb \in R^\times$ are units. But then
$\frac ra(X-t) \in R[X]$ implies $t \in R$, hence $R$ is integrally
closed.
\end{bew}

Thus, in order to transfer the Riesz property from $R$ to $R[X]$, we
need $R$ to be integrally closed. As it turns out, however, we do not
need to assume this for the passage from the lifting lemma to the
Riesz property of $R[X]$, since it follows from the lifting lemma:

\begin{prop}\label{ePGA}
Let $R$ be an integral domain. If Gauss's lifting lemma holds in
$R[X]$, then $R$ is integrally closed.
\end{prop}

\begin{bew}
Let $K$ be the quotient field of $R$ and let $t \in K$ be integral
over $R$, i.e., a root of a polynomial
$f(X) = X^n + a_{n-1}X^{n-1} + \ldots + a_0 \in R[X]$. Then there is
a factorization $f(X) = (X-t)g(X)$ in $K[X]$ with
$X-t, g(X) \in K[X]$. By Gauss's lifting lemma there is a
$q \in K^\times$ such that $q(X-t)$ and $\frac1q g(X)$ both lie in
$R[X]$. Since both $X-t$ and $g(X)$ are monic, i.e., have leading
coefficient $1$, we must have $q \in R$ and $\frac1q \in R$. Thus
$q \in R^\times$ is a unit, and this implies $X-t, g(X) \in R[X]$.
Hence $t \in R$, and $R$ is integrally closed.
\end{bew}

In particular, Gauss's lifting lemma does not hold over
$R = \Z[\sqrt{-3}\,]$: the factorization
$X^2 + X + 1 = (X-\rho)(X-\rho^2)$ with
$\rho = \frac{-1+\sqrt{-3}}2 \in K = \Q(\sqrt{-3}\,)$ clearly cannot
be lifted to a factorization in $R[X]$.

\section{Cohn's permanence theorem}\label{eAb12}

Now back to our problem: we want to show that the validity of Gauss's
lifting lemma in $R[X]$ implies that $R[X]$ is a Riesz domain: we
have already seen that $R$ is then a Riesz domain and integrally
closed; from this we would like to conclude that $R[X]$ is a Riesz
domain as well. The first step in this direction is the following

\begin{prop}\label{ePprimRX}
Let $R$ be an integrally closed Riesz domain. Then all $a \in R$ are
also primal in $R[X]$.
\end{prop}

\begin{bew}
Let $a \in R \setminus \{0\}$ and $a \mid b(X)c(X)$ for polynomials
$b, c \in R[X] \setminus \{0\}$. By Dedekind's Prague theorem this
implies\footnote{For elements $a, b, c \in R$, the relation
$a \mid bc$ in $R$ is equivalent to the one in $R[X]$, so that we
need not distinguish between $a \mid_R bc$ and
$a \mid_{R[X]} bc$.} $a \mid b_ic_j$ for all $i, j$. Since $R$ is a
Riesz domain, its divisibility group $G(R)$ is a Riesz group, and by
the Euclidean refinement property from Part I (applied to the nonzero
coefficients; coefficients $0$ are divisible by every element) there
is a factorization $a = rs$ with $r \mid b_i$ and $s \mid c_j$ for
all $0 \le i \le m$ and $0 \le j \le n$. Hence $r \mid b(X)$ and
$s \mid c(X)$; in other words: $a$ is also primal in $R[X]$.
\end{bew}

\noindent{\bf Remark.} Here it does not suffice to assume only that
the single element $a$ is primal in $R$: the Euclidean refinement
property $E_{m,n}$ uses the interpolation property of the whole group
$G(R)$, not just the primality of $a$.
\medskip

From our discussion above, the following theorem emerges:

\begin{thm}\label{eSRX}
Let $R$ be an integral domain. Then the following properties are
equivalent:
\begin{enumerate}
\item $R[X]$ is a Riesz domain.
\item Gauss's lifting lemma holds in $R[X]$.
\item $R$ is an integrally closed Riesz domain.
\item $R[X]$ is an integrally closed Riesz domain.
\end{enumerate}
\end{thm}

Integrally closed Riesz domains are also called (following Cohn)
Schreier domains.

So far we have proved (1) $\impl$ (2) $\impl$ (3), and of course
(4) $\impl$ (1) is trivial. What remains to be proved is the
``permanence theorem'' (3) $\impl$ (4) for the Schreier property.
Since the other implications are fairly straightforward, the fact
that the lifting lemma holds exactly in Schreier domains (due to
McAdam \& Rush \cite{MR}) is thus essentially equivalent to Cohn's
permanence theorem (3) $\impl$ (4) (cf.\ \cite[Thm. 9.12]{MM}).

There are some properties that are inherited by the polynomial ring
$R[X]$ from the ring $R$. It is well known, for example, that $R[X]$
is factorial whenever $R$ is\footnote{This theorem is often called
Gauss's theorem, although Gauss neither formulated nor proved it. It
was probably first explicitly formulated and proved by Hensel
\cite{Hen27}, although special cases such as the UFD property of the
polynomial rings $\Z[X_1,\ldots,X_n]$ were of course already known to
Kronecker (the theorem may also have appeared in algebra lectures,
e.g.\ by Emil Artin and Emmy Noether; in any case it can be found a
little later in van der Waerden's book \cite{vdWA}, which was based
on these lectures). The name ``Gauss's theorem'' presumably stems
from van der Waerden's remark in \cite{vdWA} that the proof to be
presented goes back to Gauss. Indeed, he proves the theorem with the
help of Gauss's Lemma.}. The property of being a GCD domain (or of
being integrally closed) also carries over from $R$ to $R[X]$. That
the Riesz property does not always carry over follows from
Prop.\ \ref{ePRX} (and an example of a ring that is a Riesz domain but
not integrally closed).

While several simple proofs of ``Gauss's theorem'' are known, the
proof of (3) $\impl$ (4) seems to lie deeper; since I only know
Cohn's proof\footnote{The proof given here is structured differently
from Cohn's (which proceeds via a result analogous to Nagata's
theorem (see, e.g., \cite{LL}) establishing a connection between the
Schreier property of a ring $R$ and of its localizations $R_S$), but
essentially uses the same tools: interpolation and refinement
properties in Riesz domains, and Dedekind's Prague theorem.}
\cite{Cohn} and its version by Anderson \& Zafrullah\footnote{Like
ours, their proof uses the Euclidean refinement property.} \cite{AZ},
I warmly encourage the readers to think about further and perhaps
more direct proofs of this theorem.

In order to show that $R[X]$ is a Riesz domain under suitable
hypotheses, we have to deduce from $a \mid bc$ for polynomials
$a, b, c \in R[X]$ the existence of a factorization $a = rs$ with
$r, s \in R[X]$ such that $r \mid b$ and $s \mid c$. In GCD domains
we could simply put $r = \gcd(a,b)$; in our case this also works if
we pass from $R[X]$ to the Euclidean ring $K[X]$, where $K$ is the
quotient field of $R$. This gives a factorization $a = \rho\sigma$
with $\rho, \sigma \in K[X]$ and $\rho \mid b$, $\sigma \mid c$ in
$K[X]$. In order to obtain a factorization in $R[X]$ from this, it
is natural to apply Gauss's lifting lemma. We therefore first prove
the direction (3) $\impl$ (2) of Theorem \ref{eSRX}:

\begin{prop}
If $R$ is an integrally closed Riesz domain, then Gauss's lifting
lemma holds in $R[X]$.
\end{prop}

\begin{bew}
Let $f \in R[X]$ and $f = gh$ with $g, h \in K[X] \setminus \{0\}$.
Then there exist $b, c \in R \setminus \{0\}$ with
$bg, ch \in R[X]$, and with $a = bc$ we have $af = (bg)(ch)$, i.e.,
$a \mid (bg)(ch)$ in $R[X]$. By the Prague theorem,
$a \mid (bg)_i(ch)_j$ for all $i, j$, where $(bg)_i = bg_i$ and
$(ch)_j = ch_j$ denote the coefficients of $bg$ and $ch$. The
Euclidean refinement property in the Riesz group $G(R)$ gives us a
factorization $a = rs$ with $r \mid bg_i$ and $s \mid ch_j$ for all
$i, j$. Thus $G = \frac br g$ and $H = \frac cs h$ are polynomials in
$R[X]$, and $GH = \frac{bc}{rs} gh = f$. This is what we wanted to
show.
\end{bew}

We now sharpen the statement of Lemma \ref{eLeGLR} and prove
(2) $\impl$ (1) of Theorem \ref{eSRX}; since we already know that the
validity of Gauss's lifting lemma in $R[X]$ implies that $R$, and
hence also $R[X]$, is integrally closed, we have thus also proved
(2) $\impl$ (4) and hence all of Theorem \ref{eSRX}.

\begin{prop}
If Gauss's lifting lemma holds in $R[X]$, then $R[X]$ is a Riesz
domain.
\end{prop}

\begin{bew}
By Lemma \ref{eLeGLR} and Prop.\ \ref{ePGA}, $R$ is an integrally
closed Riesz domain; in particular, Dedekind's Prague theorem holds
in $R[X]$ (Theorem \ref{eGLDP}), and the divisibility group
$G(R) = K^\times/R^\times$ is a Riesz group, hence has the
interpolation property $R_{m,n}$ for all $m, n$ (Part I). As before
we write $x \le y$ for $x, y \in K^\times$ if $\frac yx \in R$; in
particular, $1 \le x$ simply means $x \in R$.

Let $a \mid bc$, i.e., $ad = bc$ for polynomials
$a, b, c, d \in R[X]$. If one of these polynomials is $0$, the claim
is easily checked; so let $a, b, c, d \ne 0$. Since $K[X]$ is a GCD
domain, there are polynomials $r, s, t, u \in K[X]$ with $a = rs$,
$d = tu$, as well as $b = rt$ and $c = su$. Since Gauss's lifting
lemma holds in $R[X]$, we may modify $r, s, t, u$ by constant factors
in such a way that $r, s, t, u \in R[X]$ and the equations $a = rs$
and $d = tu$ still hold. The other two equations then only hold up
to constant factors: there is a $q \in K^\times$ with $b = q\,rt$,
and since $rt \cdot su = ad = bc = q\,rt \cdot c$, we then have
$c = q^{-1} su$.

Writing $q = \frac\beta\alpha$ with $\alpha, \beta \in R$, we have
$\alpha b = (\beta r) t$, i.e., $\alpha \mid (\beta r)t$ in $R[X]$,
and Dedekind's Prague theorem gives $\alpha \mid \beta r_it_j$, i.e.,
$1 \le q r_i t_j$ for all nonzero coefficients $r_i$ of $r$ and $t_j$
of $t$. Similarly, $\beta c = (\alpha s)u$ implies
$1 \le q^{-1} s_ku_l$ for all nonzero coefficients $s_k$ of $s$ and
$u_l$ of $u$.

We now look for a $\kappa \in K^\times$ such that
$$ r' = \kappa r, \quad s' = \kappa^{-1} s, \quad
   t' = q\kappa^{-1} t, \quad u' = q^{-1}\kappa u $$
lie in $R[X]$. Then $a = rs = r's'$, $d = tu = t'u'$,
$b = q\,rt = r't'$ and $c = q^{-1}su = s'u'$, and $a \mid bc$ has
been decomposed in the desired form. The conditions
$r', s', t', u' \in R[X]$ read
$$ r_i^{-1} \le \kappa, \quad q u_l^{-1} \le \kappa, \quad
   \kappa \le s_k, \quad \kappa \le q t_j $$
for all $i, j, k, l$ (with nonzero coefficients). By the Riesz
interpolation property such a $\kappa$ exists provided each of the
lower bounds lies below each of the upper bounds. But this is the
case:
\begin{itemize}
\item $r_i^{-1} \le s_k$ means $1 \le r_is_k$, which holds since
      $r_is_k \in R$;
\item $r_i^{-1} \le q t_j$ means $1 \le q r_i t_j$, see above;
\item $q u_l^{-1} \le s_k$ means $1 \le q^{-1} s_k u_l$, see above;
\item $q u_l^{-1} \le q t_j$ means $1 \le t_j u_l$, which holds since
      $t_ju_l \in R$.
\end{itemize}
This completes the proof.
\end{bew}

Let us point out once more that a conceptual proof of Cohn's
permanence theorem would be very welcome.

\section{Other forms of Gauss's Lemma}\label{eAb13}

In many sources, Gauss's Lemma is formulated differently. If $R$ is
an integral domain, a polynomial $f \in R[X]$ is called
{\em primitive} if every common divisor of the coefficients of $f$
is a unit in $R$. If $R$ is a GCD domain and
$f(X) = f_0+f_1X + \ldots + f_nX^n$, then
$\ct(f) = \gcd(f_0, f_1, \ldots, f_n)$ is called the content of $f$.
As Dedekind observed (in the case $R = \Z$; see \cite{Ded}), Gauss's
lifting lemma can then be formulated as follows.

\begin{lem}\label{eLggt}
If $R$ is a factorial domain (or just a GCD domain), then the product
of primitive polynomials in $R[X]$ is again primitive. More
precisely: for $f, g \in R[X]$ we have $\ct(fg) = \ct(f) \ct(g)$ (up
to units).
\end{lem}

The content $\ct(f)$ of polynomials has the property
$\ct(ag) = a \ct(g)$ for $a \in R$ and $g \in R[X]$; in particular,
$a \mid f$ for $f \in R[X]$ always implies $a \mid \ct(f)$.

The lifting lemma for GCD domains now follows like this: if $F = gh$
for polynomials $F \in R[X]$ and $g, h \in K[X]$, there exist
$a, b \in R \setminus \{0\}$ with $ag, bh \in R[X]$. Thus
$abF = (ag)(bh)$ is a factorization in $R[X]$, and
$ab\ct(F) = \ct(ag) \ct(bh)$. With $r = \gcd(ab,\ct(ag))$ and
$s = \frac{ab}r$ we then have $ab = rs$ with $r \mid ag$; since
$s = \frac{ab}r$ and $\frac{\ct(ag)}r$ are coprime and
$s \mid \frac{\ct(ag)}r \ct(bh)$, we also get $s \mid \ct(bh)$, i.e.,
$s \mid bh$. Thus $G = \frac ar g$ and $H = \frac bs h = \frac ra h$
are polynomials in $R[X]$ with $F=GH$.

\medskip\noindent{\bf Definition.}
We call $R$ a $P^4$-domain (for $P^4 = PPPP$: products of primitive
polynomials are primitive) if the product of two primitive
polynomials in $R[X]$ is again primitive.
\medskip

We have seen that Gauss's lifting lemma for monic polynomials holds
in $R[X]$ if and only if $R$ (and hence $R[X]$) is integrally closed,
or, equivalently, if Dedekind's Prague theorem holds. Similarly,
Gauss's lifting lemma for arbitrary polynomials holds exactly in
Schreier domains. It is also not hard to show that $P^4$-domains lie
between Riesz domains and AP domains (see Fig.\ \ref{eFR}; in fact, AP
already follows from the validity of the $P^4$ property for linear
polynomials). Finally, Gauss's Lemma in the form of
Lemma \ref{eLggt} holds if and only if $R$ is a GCD domain.

Now in arbitrary integral domains one can consider the ideal
$\fc(f)$ generated by the coefficients of $f$ and ask under which
conditions the equation $\fc(fg) = \fc(f)\fc(g)$ holds. The answer is
given by the following

\begin{thm}
For an integral domain $R$ the following statements are equivalent:
\begin{enumerate}
\item Every finitely generated ideal is projective.
\item Every finitely generated ideal $\ne (0)$ is invertible.
\item For every prime ideal $P \ne (0)$, the localization $R_P$ is a
      valuation ring.
\item Every ring $S$ with $R \subseteq S \subseteq K$, where $K$ is
      the quotient field of $R$, is integrally closed.
\item $\fc(fg) = \fc(f)\fc(g)$ for all polynomials $f, g \in R[X]$.
\end{enumerate}
\end{thm}

\begin{bew}
See \cite[Ch. IV]{Gil}.
\end{bew}

Thus this version of Gauss's Lemma holds if and only if the ring of
coefficients is a Pr\"ufer domain. Pr\"ufer himself (see
\cite[p. 23]{Pru}) called the theorem $\fc(fg) = \fc(f)\fc(g)$ the
``theorem of Gauss--Kronecker--Dedekind''.

\section*{Last but not least \ldots}
The motivation for writing this article essentially stems from my
reading of \cite{MM} (which deals with Gauss's lifting lemma in
Schreier domains) and \cite{PR} (where the gap in Euclid's proof of
his Prop.\ VII.19 is discussed) and from the wish to present what is
said there in a somewhat different way. I would also like to thank
Bill Dubuque for various postings in sci.math, which were very
helpful to me, especially at the beginning. Finally, I owe numerous
useful and helpful remarks to Heinz L\"uneburg and Klaus Volkert, as
well as to the two anonymous referees.

\bigskip

This article was published in Math. Semesterberichte  \textbf{56}
Bd. 56 (2009), S. 39--51. The English translation was done with the help of
by Claude (anthropic), which also closed the gap in the proof of
Theorem~\ref{eGLDP}.

\end{otherlanguage*}

\end{document}